\documentclass{amsart}
\usepackage{amssymb,latexsym}

\newtheorem{definition}{Definition}
\newtheorem{theorem}{Theorem}
\newtheorem{lemma}{Lemma}

\newtheorem{proposition}{Proposition}

\newcommand{\G}{{\mathcal{G}_\infty}}
\newcommand{\Li}{{\mathcal{L}_\infty}}
\newcommand{\A}{{\mathcal{A}_\infty}}
\newcommand{\Bi}{{\mathcal{B}_\infty}}
\newcommand{\bi}{{\frak b}}
\newcommand{\di}{{\frak d}}
\newcommand{\K}{{\bf k}}
\newcommand{\N}{{\mathbb N}}
\newcommand{\Z}{{\mathbb Z}}
\newcommand{\B}{{\frak B}}
\newcommand{\s}{{s}}
\newcommand{\2}{{s^{-1}}}
\newcommand{\Id}{{\texttt{Id}}}

\newcommand{\Deg}{{\texttt{deg}}}

\newcommand{\HC}{{\texttt{C}^*}}
\newcommand{\e}{{\texttt{Exp}}}
\newcommand{\Hom}{{\texttt{Hom}_\K}}
\newcommand{\tensor}{{\otimes_\K}}

\newcommand{\HOM}{{\underline{\texttt{Hom}}}}

\newcommand{\F}{{\texttt{F}}}

\newcommand{\T}{{\texttt{T}}}
\newcommand{\m}{{\texttt{m}}}
\newcommand{\zh}{{g}}
\newcommand{\Zh}{{G}}
\newcommand{\sh}{{h}}
\newcommand{\Sh}{{H}}
\newcommand{\g}{{\frak{\zh}}}
\newcommand{\h}{{\frak{\sh}}}
\newcommand{\Ha}{{\frak{H}}}
\newcommand{\al}{{\alpha}}
\newcommand{\be}{{\beta}}
\newcommand{\ga}{{\gamma}}
\newcommand{\Ga}{{\Gamma}}
\newcommand{\de}{{\delta}}

\newcommand{\dih}{{\widehat{\di}}}
\newcommand{\Sm}{{\underset{m>0}{\Sigma}}}
\newcommand{\SM}{{\underset{m>0}{\prod}}}
\newcommand{\Om}{{\underset{m>0}{\bigoplus}}}

\newcommand{\tensorM}{{^{\tensor^m}}}

\newcommand{\RI}{{\m_R}}
\newcommand{\R}{{R}}

\newcommand{\Aone}{{A}}
\newcommand{\Atwo}{{B}}
\newcommand{\Btwo}{{B'}}
\newcommand{\AD}{{D}}
\newcommand{\ADD}{{\AD'}}
\newcommand{\aone}{{a}}
\newcommand{\atwo}{{b}}

\newcommand{\TG}{{\T^c(\GOneTwo)}}
\newcommand{\THH}{{\T^c(\HOneTwo)}}

\newcommand{\gOneTwo}{{\psi}}

\newcommand{\gOne}{{\zh}}
\newcommand{\gTwo}{{\sh}}
\newcommand{\gOneHat}{{\widehat{\gOne}}}
\newcommand{\gTwoHat}{{\widehat{\gTwo}}}
\newcommand{\gOneTwoTilde}{{\widetilde{\gOneTwo}}}
\newcommand{\gap}{{\ga_p}}
\newcommand{\gaq}{{\ga_q}}

\newcommand{\gam}{{G_m}}
\newcommand{\gao}{{G_1}}
\newcommand{\gat}{{G_2}}
\newcommand{\gath}{{G_3}}
\newcommand{\go}{{\GMe_0}}
\newcommand{\goo}{{\zh_0}}
\newcommand{\got}{{\sh_0}}
\newcommand{\goot}{{\psi_0}}
\newcommand{\gaz}{{\ga_0}}
\newcommand{\Dg}{{d^\go}}

\newcommand{\GOneTwo}{{\Psi}}
\newcommand{\GOne}{{\g}}
\newcommand{\GTwo}{{\h}}
\newcommand{\HTwo}{{\h'}}

\newcommand{\HOneTwo}{{\Psi'}}
\newcommand{\GG}{{\Ga}}
\newcommand{\GGP}{{\GG_p}}
\newcommand{\GMorphism}{{\frak{L}(\Aone,\Atwo)}}
\newcommand{\GMe}{{\frak{l}}}
\newcommand{\GDefMor}{{\frak{L}(\mor)}}
\newcommand{\BDefMor}{{\Bi(\mor)}}

\newcommand{\BDefMorr}{{\Bi(\morr)}}

\newcommand{\BGG}{{\overline{\GG}}}
\newcommand{\BHH}{{\BGG'}}
\newcommand{\BGGP}{{\overline{\GG}_p}}

\newcommand{\BGGone}{{\overline{\GG}_1}}
\newcommand{\BGGQQ}{{\overline{\GG}_{2q}}}

\newcommand{\mor}{{f}}

\newcommand{\morr}{{f'}}

\newcommand{\Del}{{\texttt{Del}}}
\newcommand{\DEL}{{\underline{\texttt{Del}}}}
\newcommand{\DEFf}{{\underline{\texttt{Def}}(\mor)}}
\newcommand{\Deff}{{\texttt{Def}(\mor)}}
\newcommand{\Delg}{{\texttt{Del}(\GDefMor)}}
\newcommand{\DELg}{{\underline{\texttt{Del}}(\GDefMor)}}

\newcommand{\Fr}{{\F_\R}}
\newcommand{\FR}{{\underline{\F}_\R}}

\newcommand{\dgart}{{\texttt{dgart}}}

\newcommand{\CacL}{{\be}}
\newcommand{\CacR}{{\be'}}
\newcommand{\Calg}{{C}}
\newcommand{\calg}{{c}}
\newcommand{\Dc}{{d_\Calg}}
\newcommand{\Dm}{{d_\Modu}}
\newcommand{\Da}{{d_\Aalg}}
\newcommand{\Balg}{{B}}
\newcommand{\balg}{{b}}
\newcommand{\Aalg}{{A}}
\newcommand{\aalg}{{a}}
\newcommand{\Modu}{{W}}
\newcommand{\modu}{{w}}
\newcommand{\Endo}{{\frak{E}}}
\newcommand{\yeh}{{\widehat{\ye}}}
\newcommand{\ye}{{\frak{e}}}
\newcommand{\mo}{{\mu}}

\newcommand{\TC}{{\T^c}}
\newcommand{\TA}{{\T}}

\usepackage{thumbpdf}
\usepackage{hyperref}

\begin{document}
\title{$\G$-structure on deformation complex of a morphism}
\author{Dennis V. Borisov}
\date{\today}
\maketitle

\begin{abstract}
$\G$-structure is shown to exist on the deformation complex of a morphism of associative algebras. The main step of the construction is
extension of a $\Bi$-algebra by an associative algebra. Actions of $\Bi$-algebras on associative and $\Bi$-algebras are analyzed, extensions of
$\Bi$-algebras by associative and $\Bi$-algebras, that they act upon, are constructed. The resulting $\G$-algebra on the deformation complex of
a morphism is shown to be quasi-isomorphic to the $\G$-algebra on deformation complex of the corresponding diagram algebra.
\end{abstract}

\tableofcontents

\section{Introduction}
In this paper we consider simultaneous deformations of a morphism of associative dg algebras, together with its domain and codomain. In
\cite{Bor} it is shown, that the Hochschild complex of a morphism, considered in \cite{GS2}, \cite{GS3}, is not a dg Lie algebra, but is a
proper $\Li$-algebra. In this paper we seek to extend this $\Li$-structure to a $\G$-structure (\cite{TT} section 1). We work over a field $\K$
of characteristic zero.

In case of deformations of a single associative algebra, one of the methods, used to construct a $\G$-algebra on the deformation complex, is
dequantization of the $\Bi$-structure (\cite{TT} section 3, \cite{Hin2} sections 6,7). We will do the same in case of deformations of a
morphism. With this method the resulting $\G$-algebra has a dg Lie algebra as its underlying $\Li$-algebra. Therefore our construction is not
really an extension of the $\Li$-algebra from \cite{Bor}, but an extension of a quasi-isomorphic dg Lie algebra.

So our goal is to construct a $\Bi$-algebra. Here are the main steps. Let $\mor:\Aone\rightarrow\Atwo$ be a morphism of dg associative algebras
we wish to deform. We denote
    $$\GOne:=\Hom(\Om(\s\Aone)\tensorM,\s\Aone),\quad\GTwo:=\Hom(\Om(\s\Atwo)\tensorM,\s\Atwo),$$
where $\s$ stands for the suspension functor on the category of dg $\K$-spaces. Here we consider $\mor$ as an $\A$-morphism, i.e. we work with
its bar construction, hence the suspension. It is well known that $\2\GOne$, $\2\GTwo$ are $\Bi$-algebras (\cite{Get} subsection 5.2).

These $\Bi$-algebras correspond to Hopf algebras, that represent automorphism groups of $\Aone$ and $\Atwo$. The two groups act on the space of
morphisms from $\Aone$ to $\Atwo$. This space is represented by the following coassociative coalgebra without counit
    $$\underset{n>0}{\bigoplus}(\GOneTwo)^{\tensor^n},\quad\text{where }\GOneTwo=\Hom(\Om(\s\Aone)\tensorM,\s\Atwo).$$
Consequently, the two $\Bi$-algebras act on this coassociative coalgebra, $\2\GOne$ - from the right, $\2\GTwo$ - from the left. The actions are
by compositions of multilinear maps.

The notion of actions of $\Bi$-algebras on coassociative coalgebras, that we use, is a translation into $\Bi$-language of the concept of module
coalgebras over the corresponding Hopf algebras (e.g. \cite{COZ} subsection 1.3). Similarly we have the notion of an associative algebra over a
$\Bi$-algebra.

Associative algebras and coassociative coalgebras are connected by the bar and the cobar constructions. Since cobar construction is a functor,
it follows that if we have an action of a $\Bi$-algebra on a coassociative coalgebra, its cobar construction also carries an action of the same
$\Bi$-algebra. Therefore from $\underset{n>0}{\bigoplus}(\GOneTwo)^{\tensor^n}$ we get a dg associative algebra on which $\2\GOne$ and $\2\GTwo$
act. We denote this associative algebra by $\GG$.

Once we have a $\Bi$-algebra $\mathcal B$ and an associative algebra $Z$ that it acts upon, we might want to extend $\mathcal B$ by $Z$, i.e. to
find a $\Bi$-algebra $\mathcal{B}'$ and two morphisms of $\Bi$-algebras (recall that desuspension of an associative algebra is a $\Bi$-algebra):
    $$\2 Z\rightarrow\mathcal{B}'\rightarrow\mathcal B,$$
s.t. the underlying morphisms of $\K$-complexes make up a short exact sequence, and operations of $\mathcal{B}'$ extend the action of $\mathcal
B$ on $Z$. We show that such an extension is always possible.

When we apply this to the action of $\2\GOne$ on $\GG$, we get a $\Bi$-structure on $\2(\GG\oplus\GOne)$. Now we have an action of $\2\GTwo$ on
this $\Bi$-algebra (extension of the action on $\2\GG$ by the trivial action on $\2\GTwo$). Extensions of $\Bi$-algebras by $\Bi$-algebras, that
they act upon, are not very different from extensions by associative algebras. So, extending $\2\GTwo$ by $\2(\GG\oplus\GOne)$, we get a
$\Bi$-structure on $\2(\GG\oplus\GOne\oplus\GTwo)$. This is almost what we wanted. We have to complete $\GG$ with respect to the filtration by
tensor powers of $\GOneTwo$. Dequantization of the completed $\Bi$-algebra is the $\G$-algebra describing deformations of $\mor$ (as usual
dequantization depends on the choice of a Drinfeld associator, \cite{EK} subsection 2.4).

In \cite{GS1} section 23, it is shown, that for every diagram of associative algebras there is one associative algebra, called the diagram
algebra, whose Hochschild cohomology is isomorphic to that of the diagram. We show that the $\Bi$-algebra, that we construct on the deformation
complex of a morphism, is quasi-isomorphic to the usual $\Bi$-algebra on the deformation complex of the corresponding diagram algebra. Since
dequantization is a functor, the $\G$-algebras are quasi-isomorphic as well.

Here is the structure of the paper. Section \ref{SectionRepresentations} presents the necessary results about $\Bi$-algebras. In subsection
\ref{SubsectionDeformations} we recall the definition of $\Bi$-algebras, and describe the well known technique of deforming $\Bi$-algebras by
elements, that are solutions of the Maurer-Cartan equation.

In subsection \ref{Calgebras} we analyze algebras and coalgebras over $\Bi$-algebras. We show that actions of $\Bi$-algebras on (co)associative
(co)algebras are equivalent to $\Bi$-morphisms into deformation complexes of these (co)algebras. We prove that cobar and bar constructions
transport actions of $\Bi$-algebras. At the end we give the definition of a $\Bi$-algebra over a $\Bi$-algebra.

In subsection \ref{Ext} we define extensions of $\Bi$-algebras by associative algebras and by $\Bi$-algebras, and prove that it is always
possible to extend, by presenting an explicit construction. Finally we show that if we have two $\Bi$-algebras acting on one associative
algebra, s.t. these actions commute, we can extend the direct product of these $\Bi$-algebras by the associative algebra, that they act on.

Section \ref{SectionMorphisms} contains the main results of the paper. In subsection \ref{SubsectionGstructure} we apply the machinery of
section \ref{SectionRepresentations} to the case of a morphism between associative dg algebras. We get a $\Bi$-algebra, and therefore a
$\G$-algebra (depending on the choice of a Drinfeld associator).

In subsection \ref{SubsectionLieStructure} we prove that if the morphism, that we start with, is between non-positively graded algebras, then
the underlying dg Lie algebra of the $\G$-algebra, that we have constructed in subsection \ref{SubsectionGstructure}, would be the correct one
to describe deformations of that morphism. As usual, the condition on grading comes from the fact that almost free $\Z$-graded algebras do not
have to be cofibrant.

In subsection \ref{SubsectionDiagram} we show that the $\Bi$-algebra from subsection \ref{SubsectionGstructure} is quasi-isomorphic to the usual
$\Bi$-algebra on the Hochschild complex of the diagram algebra of the morphism.

\textbf{Notation}
    We fix a field $\K$ of characteristic 0. For a $\K$-space $A$ we denote by $\TA(A)$ the free associative algebra, generated by $A$, and by
    $\TC(A)$ the cofree coassociative coalgebra, cogenerated by $A$ (cofree in the category of coalgebras, that are cocomplete with respect to
    the filtration by primitives, see e.g. \cite{LM} page 2150), i.e. $\TC(A):=\underset{n>0}{\bigoplus}A^{\tensor^n}$.

    Working with Hochschild cochains $\HC(A,A)$, we denote by $\widehat{\alpha}$ the coderivation on $\TC(\s A)$ generated by $\alpha\in\HC(A,A)$,
    and by $\widetilde{\gamma}$ (for $\gamma\in\HC(A,B)$) the morphism of coalgebras $\TC(\s A)\rightarrow \TC(\s B)$.

    For typographical reasons, for a multi-linear map $\alpha$,
    instead of $\alpha(a_1\tensor...\tensor a_n)$ we write
    $\alpha(a_1,..,a_n)$. For two (or more) maps $\alpha_1$,
    $\alpha_2$ we denote by $\alpha_1\tensor\alpha_2$ the map,
    that can take value on $a_1\tensor..\tensor a_n$ for all
    $n>1$, i.e.
        $$\alpha_1\tensor\alpha_2(a_1,..,a_n):=\underset{i<j}{\Sigma}\pm
        a_1\tensor..\alpha_1(a_i)..\alpha_2(a_j)..\tensor a_n,$$
    with the signs given by the Koszul sign rule.

    Differentials raise degree. When we write a homogeneous element of a module as an exponent, we mean its parity.

\section{$\Bi$-algebras and their representations}\label{SectionRepresentations}
\subsection{Deformations of $\Bi$-algebras}\label{SubsectionDeformations}$\quad$\\
In this subsection we describe deformations of $\Bi$-algebras, given by elements that satisfy a $\Bi$-analog of the Maurer-Cartan equation.

First we recall the definition of $\Bi$-algebras. For a graded $\K$-space $\Balg$, a $\Bi$-structure on it is a dg Hopf structure on the cofree
coalgebra, cogenerated by $\s\Balg$. In order to have more transparent formulas we will work with $\Bi$-structures on $\2\Balg$.
\begin{definition}\label{Balgebra}(\cite{Get} subsection 5.2, \cite{V} subsection 2.2)
    Let $\2\Balg$ be a graded $\K$-space. A \underline{$\Bi$-structure on $\2\Balg$} is a set of $\K$-linear maps
        $$\{\bi_{m,n}:\Balg^{\tensor^{m+n}}\rightarrow \Balg\}_{m,n\geq 0},\quad\{\di_m:\Balg\tensorM\rightarrow\Balg\}_{m>0},$$
        $$\Deg(\bi_{m,n})=0,\text{ }\Deg(\di_m)=1,\text{ }\bi_{0,1}=\bi_{1,0}=\Id_{\Balg},\text{ }\bi_{0,m}=\bi_{m,0}=0\text{ for }m\neq 1,$$
    such that $\{\di_m\}_{m>0}$ is an $\A$-structure on $\2\Balg$, and if we denote
        $$\bi_p( u_1,..,u_m;v_1,..,v_n):=\underset{T}{\Sigma}(-1)^{\sigma(i,j)}
        \underset{1\leq q\leq p}{\bigotimes}\bi_{i_q,j_q}(u_{\overline{i}_{q-1}+1},..,u_{\overline{i}_q},
        v_{\overline{j}_{q-1}+1},..,v_{\overline{j}_q})$$
    (where $u_i,v_j\in\Balg$, $\overline{i}_q=\underset{1\leq r\leq q}{\Sigma}i_r$, $\overline{j}_q=\underset{1\leq r\leq q}{\Sigma}j_r$,
    $T=\{i_q,j_q\in\Z_{\geq 0}$, $\underset{1\leq q\leq p}{\Sigma}i_q=m$, $\underset{1\leq q\leq p}{\Sigma}j_q=n\}$,
    and $\sigma(i,j)$ is the sign of the permutation from $u_1,..,u_m,v_1,..,v_n$ to
    $u_1,..,u_{i_1},v_1,..,v_{j_1},..,u_{m-i_p+1},..,u_m,v_{j_{n-j_p+1}},..,v_n$), then the following would hold for all $l,m,n\in\N$
    and $u_q,v_q,w_q\in\Balg$
    \begin{equation}\label{Associativity}
        \underset{p>0}{\Sigma}\bi_{p,n}(\bi_p( u_1,..,u_l;v_1,..,v_m),w_1,..,w_n)=\qquad\qquad\qquad\qquad\qquad
    \end{equation}
        $$\qquad\qquad\qquad\qquad\qquad=\underset{p>0}{\Sigma}\bi_{l,p}(u_1,..,u_l,\bi_p(v_1,..,v_m;w_1,..,w_n)),$$
    \begin{equation}\label{Leibnitz}
        \underset{p>0}{\Sigma}\di_p(\bi_p(u_1,..,u_m;v_1,..,v_n))=\underset{p>0}{\Sigma}(\bi_{m-p+1,n}(\dih_p(u_1,..,u_m),v_1,..,v_n)+
    \end{equation}
        $$+(-1)^{\sigma}\bi_{m,n-p+1}(u_1,..,u_m,\dih_p(v_1,..,v_n))),$$
    where $\sigma=\underset{1\leq q\leq m}{\Sigma}u_q$, and $\dih_p$ is the coderivation on $\Om\Balg\tensorM$, cogenerated by $\di_p$, i.e.
        $$\dih_p(v_1,..,v_n)=\underset{0\leq i\leq n-p}{\Sigma}(-1)^{v_1+..+v_i}v_1...\di_p(v_{i+1},..,v_{i+p})...v_n.$$
\end{definition}

It is common, when dealing with dg Lie algebras, to deform the differential by bracket with an element of degree 1, that satisfies Maurer-Cartan
equation. The following definition and lemma (\cite{GV} subsection 3.1) describe a similar technique, applied to $\Bi$-algebras.

\begin{definition}\label{Bdeformation}
    Let $(\2\Balg,\{\bi_{m,n}\},\{\di_m\})$ be a $\Bi$-algebra. Let $\balg_0\in\Balg$ of degree $1$. A \underline{deformation of
    $\{\di_m\}$} by $\balg_0$ is $\{\di^{\balg_0}_m\}_{m>0}$, where
        $$\di^{\balg_0}_m(\balg_1,..,\balg_m):=\di_m(\balg_1,..,\balg_m)+\bi_{1,m}(\balg_0,\balg_1,..,\balg_m)+
        (-1)^{\sigma}\bi_{m,1}(\balg_1,..,\balg_m,\balg_0),$$
    where $\balg_i\in\Balg$, $\sigma=\balg_1+...+\balg_m+1$.
\end{definition}
\begin{lemma}\label{BDeformation}
    Let $(\2\Balg,\{\bi_{m,n}\},\{\di_m\})$ be a $\Bi$-algebra. Let $\balg_0\in\Balg$ of degree $1$. Suppose that
            $$\di_1(\balg_0)+\bi_{1,1}(\balg_0,\balg_0)=0.$$
    Then the deformation $(\2\Balg,\{\bi_{m,n}\},\{\di^{\balg_0}_m\})$ (definition \ref{Bdeformation}) is a $\Bi$-algebra.
\end{lemma}
\textbf{Proof:}
    We have to prove that the new operations satisfy Leibnitz property (equations (\ref{Leibnitz})), and that $\{\di^{\balg_0}_m\}$ is an $\A$-structure
    on $\2\Balg$. First we show that equations (\ref{Leibnitz}) hold for $\{\bi_{m,n}\}$ and $\{\di^{\balg_0}_m\}$.

    Using equations (\ref{Associativity}) we find, that for all $u_i$, $v_i\in\Balg$, $m$, $n\in\N$
    \begin{equation}\label{1}
        \underset{p>0}{\Sigma}\di^{\balg_0}_p(\bi_p(u_1,..,u_m;v_1,..,v_n))=\underset{p>0}{\Sigma}\di_p(\bi_p(u_1,..,u_m;v_1,..,v_n))+
    \end{equation}
        $$+\underset{p>0}{\Sigma}(\bi_{p,n}(\bi_p(\balg_0;u_1,..,u_m),v_1,..,v_n)+(-1)^\sigma\bi_{m,p}(u_1,..,u_m,\bi_p(v_1,..,v_n;\balg_0))),$$
    where $\sigma=u_1+...+u_m+v_1+...+v_n+1$. By definition of $\bi_p$ (definition \ref{Balgebra}) we have
        $$\bi_p(\balg_0;u_1,..,u_m)=\underset{i=0}{\overset{p-1}{\Sigma}}(-1)^{\sigma(i)}u_1...\bi_{1,m-p+1}(\balg_0,u_{i+1},..,u_{i+m-p+1})...u_m,$$
    where $\sigma(i)=u_1+...+u_i$, and $\bi_{1,m-p+1}(\balg_0,u_{i+1},..,u_{i+m-p+1})$ is inserted between $u_i$ and $u_{i+m-p+2}$. Similarly
        $$\bi_p(u_1,..,u_m;\balg_0)=\underset{i=0}{\overset{p-1}{\Sigma}}(-1)^{\sigma'(i)}u_1...\bi_{m-p+1,1}(u_{i+1},..,u_{i+m-p+1},\balg_0)...u_m,$$
    where $\sigma'(i)=u_{i+m-p+2}+...+u_m$. Therefore, since $\sigma(i)-\sigma'(i)=u_{i+1}+...+u_{i+m-p+1}-(u_1+...+u_m)$, we have
    \begin{equation}\label{2}
        \underset{p>0}{\Sigma}(\dih_p(u_1,..,u_m)+\bi_{m-p+1}(\balg_0;u_1,..,u_m))=
    \end{equation}
        $$=\underset{p>0}{\Sigma}(\dih^{\balg_0}_p(u_1,..,u_m)+(-1)^{u_1+...+u_m}\bi_{m-p+1}(u_1,..,u_m;\balg_0)).$$
    Using equations (\ref{Associativity}) and (\ref{2}), we have
    \begin{equation}\label{3}
        \underset{p>0}{\Sigma}\bi_{p,n}(\bi_p(u_1,..,u_m;\balg_0),v_1,..,v_n)=
        \underset{p>0}{\Sigma}\bi_{m,p}(u_1,..,u_m,\bi_p(\balg_0;v_1,..,v_n))=
    \end{equation}
        $$=\underset{p>0}{\Sigma}\bi_{m,n-p+1}(u_1,..,u_m,\dih^{\balg_0}_p(v_1,..,v_n)+(-1)^{\sigma}\bi_{n-p+1}(v_1,..,v_n;\balg_0)-
        \dih_p(v_1,..,v_n)),$$
    where $\sigma=v_1+...+v_n$. Now using the fact that $\{\di_m\}$ satisfy equations (\ref{Leibnitz}), we combine equations (\ref{1}),
    (\ref{2}), (\ref{3}), to get
        $$\underset{p>0}{\Sigma}\di^{\balg_0}_p(\bi_p(u_1,..,u_m;v_1,..,v_n))=$$
        $$=\underset{p>0}{\Sigma}(\bi_{m-p+1,n}(\dih^{\balg_0}_p(u_1,..,u_m),v_1,..,v_n)+
        (-1)^{\sigma}\bi_{m,n-p+1}(u_1,..,u_m,\dih^{\balg_0}_p(v_1,..,v_n))),$$
    where $\sigma=u_1+...+u_m$, i.e. $\{\di^{\balg_0}_m\}$ satisfy equations (\ref{Leibnitz}).

    It remains to show that $\{\di^{\balg_0}_m\}$ constitute an $\A$-structure on $\2\Balg$. This is the same as to say that these multi-linear maps
    are parts of a codifferential on the cofree coassociative coalgebra, cogenerated by $\Balg$, i.e. we have to prove that the following holds
    for all $n\in\N$
        $$\underset{l+m=n+1}{\Sigma}\di^{\balg_0}_l(\dih^{\balg_0}_m(v_1,..,v_n))=0.$$
    Denote $\di'_m:=\di^{\balg_0}_m-\di_m$. Then this equation can be rewritten as
    \begin{equation}\label{4}
        \underset{l+m=n+1}{\Sigma}\di_l(\dih_m(v_1,..,v_n))+\underset{l+m=n+1}{\Sigma}(\di'_l(\dih_m(v_1,..,v_n))+\di_l(\dih'_m(v_1,..,v_n)))+
    \end{equation}
        $$+\underset{l+m=n+1}{\Sigma}\di'_l(\dih'_m(v_1,..,v_n))=0.$$
    Since $\{\di_m\}$ satisfy equations (\ref{Leibnitz}) and $\Deg(\balg_0)=1$ we have
    \begin{equation}\label{a}
        \underset{l=1}{\overset{n+1}{\Sigma}}\di_l(\bi_l(\balg_0;v_1,..,v_n))+\underset{l=1}{\overset{n}{\Sigma}}\bi_{1,n-l+1}(\balg_0,\dih_l(v_1,..,v_n))=
    \end{equation}
        $$=\bi_{1,n}(\di_1(\balg_0),v_1,..,v_n),$$
    and
    \begin{equation}\label{b}
        \underset{l=1}{\overset{n+1}{\Sigma}}\di_l(\bi_l(v_1,..,v_n;\balg_0))-\underset{l=1}{\overset{n}{\Sigma}}\bi_{n-l+1,1}(\dih_l(v_1,..,v_n),\balg_0)=
    \end{equation}
        $$=(-1)^{v_1+...+v_n}\bi_{n,1}(v_1,..,v_n,\di_1(\balg_0)).$$
    From definition of $\{\bi_m\}$ and $\{\di^{\balg_0}_m\}$ we see that
    \begin{equation}\label{c}
        \bi_l(\balg_0;v_1,..,v_n)+(-1)^{v_1+..+v_n+1}\bi_l(v_1,..,v_n;\balg_0)=\dih'_{n-l+1}(v_1,..,v_n),
    \end{equation}
    for $1\leq l\leq n+1$ (we put $\dih'_0:=0$). Combining this and equations (\ref{a}),(\ref{b}) we conclude that the second summand of the
    l.h.s. of equation (\ref{4}) equals
    \begin{equation}\label{d}
        \bi_{1,n}(\di_1(\balg_0),v_1,..,v_n)-\bi_{n,1}(v_1,..,v_n,\di_1(\balg_0)).
    \end{equation}
    From equations (\ref{Associativity}) we have
        $$(-1)^{v_1+..+v_n}\bi_{l,1}(\bi_l(\balg_0;v_1,..,v_n),\balg_0)+(-1)^{v_1+..+v_n+1}\bi_{1,l}(\balg_0,\bi_l(v_1,..,v_n;\balg_0))=0.$$
    This, equation (\ref{c}), and definition of $\{\di^{\balg_0}\}$ imply that
        $$\underset{l+m=n+1}{\Sigma}\di'_l(\dih'_m(v_1,..,v_n))=\underset{l=1}{\overset{n}{\Sigma}}(\bi_{1,l}(\balg_0,\bi_l(\balg_0;v_1,..,v_n))-
        \bi_{l,1}(\bi_l(v_1,..,v_n;\balg_0),\balg_0)).$$
    Now we can use equations (\ref{Associativity}) and the fact that $\Deg(\balg_0)=1$ to rewrite the r.h.s. of the last equation as follows
    \begin{equation}\label{e}
        \bi_{1,n}(\bi_{1,1}(\balg_0,\balg_0),v_1,..,v_n)-\bi_{n,1}(v_1,..,v_n,\bi_{1,1}(\balg_0,\balg_0)).
    \end{equation}
    So we have proved that the sum of the second and the third summands on the l.h.s. of equation (\ref{4}) equals the sum of expressions
    (\ref{d}) and (\ref{e}). This sum is zero, because $\balg_0$ is assumed to satisfy $\di_1(\balg_0)+\bi_{1,1}(\balg_0,\balg_0)=0$.

    The first summand of the l.h.s. of equation (\ref{4}) is zero because $\{\di_m\}$ is an $\A$-structure by assumption. So equation (\ref{4})
    holds and we are done.
$\blacksquare$
\subsection{Algebras and coalgebras over $\Bi$-algebras}\label{Calgebras}$\quad$\\
For any operad there is a well known notion of modules over an algebra over that operad (e.g. \cite{KM} definition 4.1). Although $\Bi$-algebras
are algebras over $\Bi$-operad, we will use a different concept of modules over them, namely the one coming from modules over Hopf algebras
(e.g. \cite{COZ} subsection 1.1). This will be very natural in our constructions, since for the $\Bi$-algebras that we will use, the
corresponding Hopf algebras are more basic objects.
\begin{definition}\label{Modules}
    Let $(\2\Balg,\{\bi_{m,n}\},\{\di_m\})$ be a $\Bi$-algebra, let $\Modu$ be a dg $\K$-space. A \underline{module structure} on $\Modu$
    over $\2\Balg$ consists of a set of $\K$-linear maps
        $$\{\be_m:\Balg^{\tensor^{m-1}}\tensor\Modu\rightarrow\Modu\}_{m>1},\quad \Deg(\be_m)=0,$$
    such that
    \begin{equation}\label{AssAction}
        \underset{p>1}{\Sigma}\be_p(\bi_{p-1}(u_1,..,u_m;v_1,..,v_n),\modu)=\be_{m+1}(u_1,..,u_m,\be_{n+1}(v_1,..,v_n,\modu)),
    \end{equation}
    \begin{equation}\label{LeibnitzAction}
        \Dm(\be_m(v_1,..,v_{m-1},\modu))=\underset{p=1}{\overset{m-1}{\Sigma}}\be_{m-p+1}(\dih_p(v_1,..,v_{m-1}),\modu)+
    \end{equation}
        $$+(-1)^{v_1+...+v_{m-1}}\be_m(v_1,..,v_{m-1},\Dm(\modu)),$$
    where $u_i,v_i\in\Balg$, $\modu\in\Modu$ and $\Dm$ is the differential on $\Modu$.
\end{definition}

Modules over $\Bi$-algebras, as in definition \ref{Modules}, can be used to define ``square-zero'' extensions of $\Bi$-algebras. We are
interested in extending $\Bi$-algebras by associative algebras. This makes sense, since desuspensions of associative algebras are
$\Bi$-algebras.

The following definition is the $\Bi$-version of module-algebras and module-coalgebras over Hopf algebras (e.g. \cite{COZ} subsection 1.3).
\begin{definition}\label{AlgebraCoalgebra}
    Let $\2\Balg$ be a $\Bi$-algebra, let $(\Aalg,\cdot,\Da)$ be a dg associative $\K$-algebra and $(\Calg,\Delta,\Dc)$ a dg coassociative $\K$-coalgebra.
    \begin{itemize}
    \item[1.] An \underline{action of $\2\Balg$ on $\Aalg$} is a set of $\K$-linear maps $\{\be_m\}_{m>1}$, s.t. they make $\Aalg$
        into a module over $\2\Balg$ (definition \ref{Modules}) and also satisfy
        \begin{equation}\label{AlgebraAction}
            \be_m(v_1,..,v_{m-1},\aalg_1\cdot\aalg_2)=
        \end{equation}
            $$=\underset{i=0}{\overset{m-1}{\Sigma}}(-1)^{\sigma(i)}\be_{i+1}(v_1,..,v_i,\aalg_1)\cdot\be_{m-i}(v_{i+1},..,v_{m-1},\aalg_2),$$
        where $\sigma(i):=\aalg_1(v_{i+1}+...+v_{m-1})$, $v_i\in\Balg$, $\aalg_i\in\Aalg$, and we set $\be_1:=\Id_\Aalg$.
    \item[2.] An \underline{action of $\2\Balg$ on $\Calg$} is a set of $\K$-linear maps $\{\be_m\}_{m>1}$, s.t. they make $\Calg$ into a
        module over $\2\Balg$ and also satisfy
        \begin{equation}\label{CoalgebraAction}
            \Delta(\be_m(v_1,..,v_{m-1},\calg))=
        \end{equation}
            $$=\Sigma\underset{i=0}{\overset{m-1}{\Sigma}}(-1)^{\sigma(i)}\be_{i+1}(v_1,..,v_i,
            \calg_1)\tensor\be_{m-i}(v_{i+1},..,v_{m-1},\calg_2),$$
        where $\calg\in\Calg$, $\Delta(\calg)=\Sigma\calg_1\tensor\calg_2$, $\sigma(i)=\calg_1(v_{i+1}+...+v_{m-1})$ and we set
        $\be_1:=\Id_\Calg$.
    \end{itemize}
\end{definition}

We would like to have a presentation of modules, algebras and coalgebras over a $\Bi$-algebra $\2\Balg$, as certain morphisms from $\2\Balg$
into $\Bi$-algebras, that are naturally associated to the modules, algebras and coalgebras. For that purpose we need the notion of a
$\Bi$-morphism.
\begin{definition}
    Let $(\2\Balg,\{\bi_{m,n}\},\{\di_m\})$, $(\2\Balg',\{\bi'_{m,n}\},\{\di'_m\})$ be two $\Bi$-algebras. A \underline{$\Bi$-morphism}
    $\2\Balg\rightarrow\2\Balg'$ consists of $\K$-linear maps of degree $0$
        $$\mo_m:\Balg\tensorM\rightarrow\Balg',\quad m>0,$$
    s.t. they generate a morphism of dg Hopf algebras $\Om\Balg\tensorM\rightarrow\Om\Balg'\tensorM$.
\end{definition}

Let $\Modu$ be a graded $\K$-space. There are two well known $\Bi$-algebras (\cite{Get}, subsection 5.2), associated to $\Modu$:
    $$\2\Endo(\Modu):=\2\Hom(\Om\Modu\tensorM,\Modu),\quad\2\Endo'(\Modu):=\2\Hom(\Modu,\Om\Modu\tensorM).$$
We will call the first one the \underline{endomorphism $\Bi$-algebra} and the second one the \underline{co-endomorphism
$\Bi$}\underline{-algebra} corresponding to $\Modu$. In both cases the $\{\di_m\}$ operations are trivial, and the $\{\bi_{m,n}\}$ operations
are given by all possible compositions of maps.

Elements of $\2\Endo(\Modu)$ correspond to coderivations of the cofree coassociative coalgebra $\TC(\Modu)$, cogenerated by $\Modu$, while
elements of $\2\Endo'(\Modu)$ correspond to derivations of the free algebra $\TA(\Modu)$, generated by $\Modu$.

With these notation $\TC(\Modu)$ is a coalgebra over $\2\Endo(\Modu)$, with the action inductively defined as follows
\begin{equation}\label{EndoAction}
    \be_2(\ye,-):=\yeh,\quad\be_n(\ye_1,..,\ye_{n-1},\modu)=0,\text{ if }n>2,\text{ }\modu\in\Modu,
\end{equation}
    $$\Delta\circ\be_m(\ye_1,..,\ye_{m-1},-)=(\underset{i=0}{\overset{m-1}{\Sigma}}\be_{i+1}(\ye_1,..,\ye_i,-)\tensor\
    \be_{m-i}(\ye_{i+1},..,\ye_{m-1},-))\circ\Delta,$$
where $\Delta$ is the comultiplication on $\TC(\Modu)$, and $\yeh$ is the coderivation on $\TC(\Modu)$, cogenerated by $\ye\in\Endo(\Modu)$.

Similarly, $\TA(\Modu)$ is an algebra over $\2\Endo'(\Modu)$, with action inductively defined as follows
\begin{equation}\label{CoEndoAction}
    \be'_2(\ye',-):=\yeh',\quad\be'_n(\ye'_1,..,\ye'_{n-1},\modu)=0,\text{ if }n>2,\text{ }\modu\in\Modu,
\end{equation}
    $$\be'_m(\ye'_1,..,\ye'_{m-1},\aalg_1\cdot\aalg_2):=\underset{i=0}{\overset{m-1}{\Sigma}}(-1)^{\sigma(i)}\be'_{i+1}(\ye'_1,..,\ye'_i,\aalg_1)
    \cdot\be'_{m-i}(\ye'_{i+1},..,\ye'_{m-1},\aalg_2),$$
where $\aalg_i\in\TA(\Modu)$, $\sigma(i):=\aalg_1(\ye'_{i+1}+..+\ye'_{m-1})$, and $\yeh'$ stands for the derivation on $\TA(\Modu)$, generated
by $\ye'\in\Endo'(\Modu)$.

In both cases equations (\ref{AssAction}) are easily checked by direct calculation, and equations (\ref{AlgebraAction}), (\ref{CoalgebraAction})
are evident from the inductive definition.

Using lemma \ref{BDeformation} we can extend these actions to almost (co)free algebras and coalgebras (co)generated by $\Modu$. That is, let
$\yeh_0$ be a codifferential on $\TC(\Modu)$, cogenerated by $\ye_0\in\Endo(\Modu)$, and let $(\2\Endo(\Modu),\{\di^{\ye_0}_m\})$ be the
deformation of $\2\Endo(\Modu)$, defined by $\ye_0$ (lemma \ref{BDeformation}). Then the same operations $\{\be_m\}$ from equations
(\ref{EndoAction}), define an action of $(\2\Endo(\Modu),\{\di_m^{\ye_0}\})$ on $(\TC(\Modu),\yeh_0)$ (definition \ref{AlgebraCoalgebra}).

Similarly, let $\yeh'_0$ be a differential on $\TA(\Modu)$, generated by $\ye'_0\in\Endo'(\Modu)$, and $(\2\Endo'(\Modu),\{\di^{\ye'_0}_m\})$
the corresponding deformation of $\2\Endo(\Modu)$. Then operations $\{\be'_m\}$ from equations (\ref{CoEndoAction}), define an action of
$(\2\Endo'(\Modu),\{\di_m^{\ye'_0}\})$ on $(\TA(\Modu),\yeh'_0)$ (definition \ref{AlgebraCoalgebra}).

In the following lemma we use $\2\Endo(\Modu)$ and $\2\Endo'(\Modu)$ to represent actions of $\Bi$-algebras on almost (co)free (co)algebras,
(co)generated by $\Modu$, as $\Bi$-morphisms.
\begin{lemma}\label{Representation}
    Let $(\TC(\Modu),\Dc)$, $(\TA(\Modu),\Da)$ be almost (co)free coalgebra and algebra, (co)generated by a graded $\K$-space $\Modu$. Let
    $(\2\Endo(\Modu),\{\di_m\})$, $(\2\Endo'(\Modu),\{\di'_m\})$ be the corresponding deformations of the endomorphism and co-endomorphism
    $\Bi$-algebras. Let $\2\Balg$ be a $\Bi$-algebra. There is a bijection between the set of actions of $\2\Balg$ on $(\TC(\Modu),\Dc)$ (definition
    \ref{AlgebraCoalgebra}) and the set of $\Bi$-morphisms $\2\Balg\rightarrow(\2\Endo(\Modu),\{\di_m\})$. Similarly, there is a bijection between
    the set of actions of $\2\Balg$ on $(\TA(\Modu),\Da)$ and the set of $\Bi$-morphisms $\2\Balg\rightarrow(\2\Endo'(\Modu),\{\di'_m\})$.
\end{lemma}
\textbf{Proof:}
    We will prove only the coalgebra case. The other case is done similarly.
    Coalgebras over a $\Bi$-algebra were defined as the $\Bi$-version of module-coalgebras over the corresponding Hopf algebra. Therefore, if we
    have an action of $(\2\Endo(\Modu)$, $\{\di_m\})$ on $(\TC(\Modu),\Dc)$ and a $\Bi$-morphism $\2\Balg\rightarrow(\2\Endo(\Modu),\{\di_m\})$,
    there is an induced action of $\2\Balg$ on $(\TC(\Modu),\Dc)$. So we get a map from the set of $\Bi$-morphisms to the set of actions.

    Suppose we have an action $\{\be_m\}$ of $\2\Balg$ on $(\TC(\Modu),\Dc)$. From equations (\ref{CoalgebraAction}) we know that $\be_2(\balg,-)$
    is a coderivation on $(\TC(\Modu),\Dc)$ for all $\balg\in\Balg$, so we have a map $\Balg\rightarrow\Endo(\Modu)$. Let $\balg_1$,
    $\balg_2\in\Balg$. From equations (\ref{CoalgebraAction}) we conclude that
    $\be_3(\balg_1,\balg_2,-)-\be_2(\balg_1,-)\tensor\be_2(\balg_2,-)$ is a coderivation on $(\TC(\Modu),\{\di_m\})$, and we have a map
    $\Balg^{\tensor^2}\rightarrow\Endo(\Modu)$. Continuing in this way we get a sequence of maps $\{\mu_m:\Balg\tensorM\rightarrow\Endo(\Modu)\}_{m>0}$.
    It is easily checked by direct computation, that this is a $\Bi$-morphism from $\2\Balg$ to $(\2\Endo(\Modu),\{\di_m\})$. So we have a map
    from the set of actions to the set of $\Bi$-morphisms. It is the inverse of the map in the opposite direction, that is described in the beginning of the proof.
$\blacksquare$

We will need a way of changing coalgebras over $\Bi$-algebras to algebras over them. We will do it by means of the cobar construction. Since it
is a functor, it maps coderivations to derivations. When $\Bi$-algebra acts on a coalgebra (definition \ref{AlgebraCoalgebra}), elements act as
coderivations, pairs of elements - as a sort of codifferential operators of order 2 and so forth. The following lemma shows that all these
actions survive the cobar functor.
\begin{lemma}\label{Cobar}
    Let $(\Calg,\Delta,\Dc)$ be a dg coassociative coalgebra and let $(\Aalg,\cdot,\Da)$ be a dg associative algebra. Let $(\2\Balg$,
    $\{\bi_{m,n}\}$,$\{\di_m\})$ be a $\Bi$-algebra. Let $\{\be_m\}$ be an action of $\2\Balg$ on $\Calg$, and $\{\be'_m\}$ - an action
    of $\2\Balg$ on $\Aalg$ (definition \ref{AlgebraCoalgebra}).
    Then the cobar construction $\Omega(\Calg)$ of $\Calg$ is an algebra over $\2\Balg$, with the operations $\{\Omega(\be_m)\}$ uniquely defined by
    \begin{equation}\label{OverBar}
        \Omega(\be_m)(\balg_1,..,\balg_{m-1},\2\calg):=\2\be_m(\balg_1,..,\balg_{m-1},\calg),
    \end{equation}
        $$\2\calg\in\2\Calg\subset\Omega(\Calg),\quad\balg_i\in\Balg.$$
    Similarly the bar construction $\B(\Aalg)$ of $\Aalg$ is a coalgebra over $\2\Balg$, with the operations $\{\B(\be'_m)\}$ uniquely defined by
        $$\B(\be'_m)(\balg_1,..\balg_{m-1},\s\aalg):=\s\be'_m(\balg_1,..,\balg_{m-1},\aalg),\quad\s\aalg\in\s\Aalg\subset\B(\Aalg),\balg_i\in\Balg.$$
\end{lemma}
\textbf{Proof:}
    We will prove only the cobar construction part, the other part is done similarly.

    By definition $\{\be_m\}$ make $\Calg$ into a module over $\Balg$. Therefore $\2\Calg$ is also a module over $\Balg$ with the action defined
    as in equation (\ref{OverBar}). From this, since $\Omega(\Calg)$ is an almost free algebra, generated by $\2\Calg$, formulas (\ref{AlgebraAction})
    inductively define a unique action of $\Balg$ on $\Omega(\Calg)$, if we forget the differential on $\Omega(\Calg)$ and the $\{\di_m\}$
    operations on $\Balg$. So we only have to check that this action satisfies equations (\ref{LeibnitzAction}) with respect to the differential
    on $\Omega(\Calg)$.

    This differential consists of two parts. The first one is an extension of $\Dc$. Since we know that $\Calg$ is a module over $\Balg$,
    equations (\ref{LeibnitzAction}) are satisfied with respect to this extension of $\Dc$. Therefore it remains to show that the rest of the
    differential on $\Omega(\Calg)$ commutes with the action.

    The second part of the differential on $\Omega(\Calg)$ is generated by
        $$\de(\2\calg):=\Sigma(-1)^{\sigma}\2\calg_1\tensor\2\calg_2,$$
    where $\Delta(\calg)=\Sigma\calg_1\tensor\calg_2$, $\sigma:=\calg_1+1$. We will show that $\de$ and the action commute on generators of
    $\Omega(\Calg)$, extension to all of $\Omega(\Calg)$ is straightforward.

    So we have to show that
    \begin{equation}\label{FD2}
        \de(\Omega(\be_m)(\balg_1,..,\balg_{m-1},\2\calg))=(-1)^{\sigma}\Omega(\be_m)(\balg_1,..,\balg_{m-1},\de(\2\calg)),
    \end{equation}
    where $\sigma=\balg_1+...+\balg_{m-1}$. Let $\Delta(\calg)=\Sigma\calg_1\tensor\calg_2$, then because $\Calg$ is a coalgebra over $\2\Balg$, we have
        $$\Delta(\be_m(\balg_1,..,\balg_{m-1},\calg))=\Sigma\underset{i=0}{\overset{m-1}{\Sigma}}(-1)^{\sigma(i)}\be_i(\balg_1,..,\balg_i,
        \calg_1)\tensor\be_{m-i}(\balg_{i+1},..,\balg_{m-1},\calg_2),$$
    where $\sigma(i)=\calg_1(\balg_{i+1}+...+\balg_{m-1})$. Therefore the left hand side of equation (\ref{FD2}) is
        $$\Sigma\underset{i=0}{\overset{m-1}{\Sigma}}(-1)^{\sigma'(i)}\2\be_i(\balg_1,..,\balg_i,\calg_1)\tensor
        \2\be_{m-i}(\balg_{i+1},..,\balg_{m-1},\calg_2),$$
    where $\sigma'(i)=\balg_1+...+\balg_i+\calg_1+1+\calg_1(\balg_{i+1}+...+\balg_{m-1})$. The right hand side of equation (\ref{FD2}) is
        $$\Sigma\underset{i=0}{\overset{m-1}{\Sigma}}(-1)^{\sigma''(i)}\Omega(\be_i)(\balg_1,..,\balg_i,\2\calg_1)\tensor
        \Omega(\be_{m-i})(\balg_{i+1},..,\balg_{m-1},\2\calg_2),$$
    where $\sigma''(i)=\calg_1+1+\balg_1+...+\balg_{m-1}+(\calg_1+1)(\balg_{i+1}+...+\balg_{m-1})$. We see that $\sigma'(i)=\sigma''(i)$
    and hence equation (\ref{FD2}) does hold.
$\blacksquare$

Together with the notion of an associative algebra over $\Bi$-algebra, we will consider $\Bi$-algebras over $\Bi$-algebras. Let
$(\2\Balg,\{\bi_{m,n}\},\{\di_m\})$, $(\2\Balg',\{\bi'_{m,n}\},\{\di'_m\})$ be $\Bi$-algebras. Let $\{\be_m\}$ be an action of $\2\Balg$ on
$\Balg'$ as a module (definition \ref{Modules}). By means of equations (\ref{CoalgebraAction}) this action extends to an action of $\2\Balg$ on
the cofree coalgebra, cogenerated by $\Balg'$. We will denote this extension by the same symbols $\{\be_m\}$.

The following definition is a particular case of a translation into $\Bi$-language of the notion of a bialgebra over a Hopf algebra.
\begin{definition}\label{BBAction}
    With the above notation, suppose that $\{\be_m\}$ make $(\TC(\Balg'),\{\di'_m\})$ into a coalgebra over $\2\Balg$
    (definition \ref{AlgebraCoalgebra}). Also assume they make $(\TC(\Balg'),\{\di'_m\})$ into an algebra over $\2\Balg$, i.e. they satisfy
        $$\be_m(\balg_1,..,\balg_{m-1},\underset{n>0}{\Sigma}\bi'_n(u_1,..,u_p;v_1,..,v_q))=$$
        $$=\underset{n>0}{\Sigma}\underset{i=0}{\overset{m-1}{\Sigma}}
            (-1)^{\sigma(i)}\bi'_n(\be_{i+1}(\balg_1,..,\balg_i,(u_1,..,u_p));\be_{m-i}(\balg_{i+1},..,\balg_{m-1},(v_1,..,v_q))),$$
    where $u_i$, $v_i\in\Balg'$, $\balg_i\in\Balg$, $\sigma(i)=(u_1+...+u_p)(\balg_{i+1}+...+\balg_{m-1})$ and we set $\be_1:=\Id_{\Balg'}$.
    In this case we will call $\2\Balg'$ a \underline{$\Bi$-algebra over a $\Bi$-algebra} $\2\Balg$.
\end{definition}

\subsection{Extensions of $\Bi$-algebras}\label{Ext}$\quad$\\
In this subsection we construct extensions of $\Bi$-algebras. First we consider extending them by associative algebras, and then by
$\Bi$-algebras.

We illustrate the procedure with an example of a particularly simple $\Bi$-algebra. Suppose we have an associative $\K$-algebra $A$, that we
consider as a $\Bi$-algebra, i.e. $\2 A$ whose only non-trivial operation is $\bi_{1,1}$. Let $M$ be an algebra over $\2 A$, s.t. the only
non-trivial component of the corresponding $\Bi$-morphism (lemma \ref{Representation}) is a $\K$-linear map $A\rightarrow\Endo'(M)$.

Now we want to define a $\Bi$-structure on $\2(M\oplus A)$, extending the action and the existing operations on $A$ and $M$. However, the action
of elements of $A$ on $M$ is by derivations, so we have
    $$a(m_1m_2)=(am_1)m_2+(-1)^{am_1}m_1(am_2),\quad a\in A,\quad m_i\in M,$$
and the second summand on the right hand side breaks associativity. Therefore we have to introduce a correction:
    $$\bi_{2,1}(m_1\tensor a\tensor m_2):=m_1(am_2).$$
The following lemma describes the same procedure in case of associative algebra over a general $\Bi$-algebra.
\begin{lemma}\label{Extensions}
    Let $(\2\Balg,\{\bi_{m,n}\},\{\di_m\})$ be a $\Bi$-algebra (definition \ref{Balgebra}), let $(\Aalg,\cdot,\Da)$ be an associative algebra and
    $\{\be_m\}$ an action of $\2\Balg$ on $\Aalg$ (definition \ref{AlgebraCoalgebra}). Define
    $\K$-linear maps
        $$\{\bi'_{m,n}:(\Aalg\oplus\Balg)^{\tensor^{m+n}}\rightarrow\Aalg\oplus\Balg\}_{m,n>0},\quad\{\di'_m:(\Aalg\oplus\Balg)^{\tensor^m}
        \rightarrow\Aalg\oplus\Balg\}_{m>0}$$
    as follows: for all $\balg_i\in\Balg$, $\aalg_i\in\Aalg$
        $$\bi'_{m,n}(\balg_1,..,\balg_{m+n}):=\bi_{m,n}(\balg_1,..,\balg_{m+n}),\quad\bi'_{1,1}(\aalg_1,\aalg_2):=\aalg_1\cdot\aalg_2,
        \quad\di'_1(\aalg):=\Da(\aalg)$$
        $$\di'_m(\balg_1,..,\balg_m):=\di_m(\balg_1,..,\balg_m),\quad
        \bi'_{m,1}(\balg_1,..,\balg_m,\aalg):=\be_{m+1}(\balg_1,..,\balg_m,\aalg),$$
        $$\bi'_{m,1}(\aalg_1,\balg_1,..,\balg_{m-1},\aalg_2):=\aalg_1\cdot\be_m(\balg_1,..,\balg_{m-1},\aalg_2),$$
    and for the rest of possibilities the values are $0$. Then if we set $\bi'_{m,0}$, $\bi'_{0,n}$ as
    in definition \ref{Balgebra}, $(\2(\Aalg\oplus\Balg),\{\bi'_{m,n}\}_{m,n\geq 0},\{\di'_m\}_{m>0})$ is a $\Bi$-algebra.
\end{lemma}
\textbf{Proof:}
    We have to prove that $\{\bi'_{m,m}\}$ satisfy associativity conditions (equations (\ref{Associativity})) and that $\{\di'_{m}\}$ have the
    Leibnitz property with respect to $\{\bi'_{m,m}\}$ (equations (\ref{Leibnitz})).

    First we show that equations (\ref{Associativity}) are satisfied. Let $(x_1...x_n)\in(\Aalg\oplus\Balg)^{\tensor^n}$. There is
    a $\bi'_{l,m}$, s.t. $l+m=n$ and $\bi'_{l,m}(x_1,..,x_n)\neq 0$ only in one of the following cases
    \begin{itemize}
    \item[1.] $n=2$, $l=m=1$, $x_1,x_2\in\Aalg$,
    \item[2.] $n\geq 2$, $l=n-1$, $m=1$, $x_1,...,x_{n-1}\in\Balg$, $x_n\in\Aalg$,
    \item[3.] $n\geq 3$, $l=n-1$, $m=1$, $x_1,x_n\in\Aalg$, $x_2,...,x_{n-1}\in\Balg$,
    \item[4.] $n\geq 2$, $l,m\geq 1$, $x_i\in\Balg$ for all $0\leq i\leq n$.
    \end{itemize}
    Equations (\ref{Associativity}) are quadratic, and to prove that they hold, we have to show that any two compositions of pairs of the above $4$
    possibilities, applied to the same argument, produce the same result. Since $\Aalg$ is an associative algebra and $\Balg$ is a
    $\Bi$-algebra, we can exclude arguments that consist solely of elements of $\Aalg$ or $\Balg$. Then, taking into account cases when
    $\{\bi'_{m,n}\}$ vanish, we are left with the following equalities to prove:
    $$\bi'_{m+1,1}(\aalg_1,\balg_1,..,\balg_m,\bi'_{1,1}(\aalg_2,\aalg_3))=\bi'_{1,1}(\bi'_{m+1,1}(\aalg_1,\balg_1,..,\balg_m,\aalg_2),\aalg_3)+$$
        $$+\bi'_{m+1,1}(\bi'_{m+1}(\aalg_1,\balg_1,..,\balg_m;\aalg_2),\aalg_3),$$
    $$\bi'_{1,1}(\aalg_1,\bi'_{m,1}(\balg_1,..,\balg_m,\aalg_2))=\bi'_{m+1,1}(\bi'_{m+1}(\aalg_1;\balg_1,..,\balg_m),\aalg_2),$$
    $$\bi'_{m,1}(\balg_1,..,\balg_m,\bi'_{1,1}(\aalg_1,\aalg_2))=\underset{p>0}{\Sigma}\bi'_{p,1}(\bi'_p(\balg_1,..,\balg_m;\aalg_1),\aalg_2),$$
    $$\bi'_{1,1}(\aalg_1,\bi'_{m+1,1}(\aalg_2,\balg_1,..,\balg_m,\aalg_3))=\bi'_{m+1,1}(\bi'_{1,1}(\aalg_1,\aalg_2),\balg_1,..,\balg_m,\aalg_3),$$
    $$\bi'_{m,1}(u_1,..,u_m,\be'_{n,1}(v_1,..,v_n,\aalg))=\underset{p>0}{\Sigma}\bi'_{p,1}(\bi'_p(u_1,..,u_m;v_1,..,v_n),\aalg),$$
    $$\bi'_{m,1}(u_1,..,u_m,\bi'_{n+1,1}(\aalg_1,v_1,..,v_n,\aalg_2))=\underset{p>0}{\Sigma}\bi'_{p,1}(\bi'_p(u_1,..,u_m;\aalg_1,v_1,..,v_n),\aalg_2),$$
    $$\bi'_{m+1,1}(\aalg_1,u_1,..,u_m,\bi'_{n,1}(v_1,..,v_n,\aalg_2))=\underset{p>0}{\Sigma}\bi'_{p,1}(\bi'_p(\aalg_1,u_1..u_m;v_1...v_n),\aalg_2),$$
    $$\bi'_{m+1,1}(\aalg_1,u_1,..,u_m,\bi'_{n+1,1}(\aalg_2,v_1,..,v_n,\aalg_3))=\underset{p>0}{\Sigma}\bi'_{p,1}(\bi'_p(\aalg_1,u_1,..,u_m;\aalg_2,
        v_1,..,v_n),\aalg_3),$$
    where $\balg_i$, $u_i$, $v_i\in\Balg$, $\aalg_i\in\Aalg$.

    These equations are easily checked by direct computation. Consider for example the last one. The left hand side of it is
        $$\bi'_{m+1,1}(\aalg_1,u_1,..,u_m,\aalg_2\cdot\be_{n+1}(v_1,..,v_n,\aalg_3))=\aalg_1\cdot\be_{m+1}(u_1,..,u_m,\aalg_2\cdot
            \be_{n+1}(v_1,..,v_n,\aalg_3))=$$
        $$=\underset{i=0}{\overset{m}{\Sigma}}(-1)^{\sigma(i)}\aalg_1\cdot\be_{i+1}(u_1,..,u_i,\aalg_2)\cdot
            \be_{m-1+1}(u_{i+1},..,u_m,\be_{n+1}(v_1,..,v_n,\aalg_3))=$$
        $$=\underset{i=0}{\overset{m}{\Sigma}}\underset{p>0}{\Sigma}(-1)^{\sigma(i)}\aalg_1\cdot\be_{i+1}(u_1,..,u_i,\aalg_2)\cdot
            \be_{p+1}(\bi_p(u_{i+1},..,u_m;v_1,..,v_n),\aalg_3),$$
    where $\sigma(i)=\aalg_2(u_{i+1}+...+u_m)$. The last expression can now be seen to coincide with the right hand side.

    So if we forget the $\di$-part, $(\2(\Aalg\oplus\Balg),\{\bi'_{m,n}\})$ is a $\Bi$-algebra. That is equivalent to the cofree
    coassociative coalgebra cogenerated by $\Aalg\oplus\Balg$ being a Hopf algebra. To finish the proof we have to show that $\{\di'_m\}$ define a
    codifferential on this Hopf algebra.

    Since $\{\di'_m\}$ is a sum of $\{\di_m\}$ on $\Balg$ and the multiplication on $\Aalg$, it is clear that its square is 0. We have to show that
    it satisfies the Leibnitz identity with respect to $\{\bi'_{m,n}\}$. Since $\Aalg$ is a module over $\Balg$, it is enough to check that
        $$\Da(\bi'_{m,1}(\aalg_1,\balg_1,..,\balg_{m-1},\aalg_2))=\bi'_{m,1}(\Da(\aalg_1),\balg_1,..,\balg_{m-1},\aalg_2))+$$
        $$+\underset{p>0}{\Sigma}(-1)^{\aalg_1}\bi'_{m-p+1,1}(\aalg_1,\dih_p(\balg_1,..,\balg_{m-1}),\aalg_2)+
            (-1)^\sigma\bi'_{m,1}(\aalg_1,\balg_1,..,\balg_{m-1},\Da(\aalg_2)),$$
    where $\sigma=\aalg_1+\balg_1+...+\balg_{m-1}$. Again, since $\Aalg$ is a dg algebra and a module over $\2\Balg$, this is true.
$\blacksquare$

As we have already mentioned, desuspensions of associative algebras are particular cases of $\Bi$-algebras, so it is natural to expect a similar
result on extending $\Bi$-algebras by general $\Bi$-algebras. The following is an explicit construction of such extensions.

Let $(\2\Balg,\{\bi_{m,n}\},\{\di_m\})$, $(\2\Balg',\{\bi'_{m,n}\},\{\di'_m\})$ be $\Bi$-algebras. Suppose $\2\Balg$ acts on $\2\Balg'$ through
$\{\be_m\}$ (definition \ref{BBAction}). Generalizing lemma \ref{Extensions}, we define $\K$-multilinear maps
$\{\bi''_{m,n}:(\Balg'\oplus\Balg)^{\tensor^{m+n}}\rightarrow\Balg'\oplus\Balg\}$ as follows: for all $\balg_i\in\Balg$, $\balg'_i\in\Balg'$
\begin{equation}\label{Baction1}
    \bi''_{m,n}(\balg_1,..,\balg_{m+n}):=\bi_{m,n}(\balg_1,..,\balg_{m+n}),\quad\di''_m(\balg_1,..,\balg_m):=\di_m(\balg_1,..,\balg_m),
\end{equation}
\begin{equation}\label{Baction2}
    \bi''_{m,n}(\balg'_1,..,\balg'_{m+n}):=\bi'_{m,n}(\balg'_1,..,\balg'_{m+n}),\quad\di''_m(\balg'_1,..,\balg'_m):=\di'_m(\balg'_1,..,\balg'_m),
\end{equation}
\begin{equation}\label{Baction3}
    \bi''_{m,n}(\balg_1,..,\balg_m,\balg'):=\be_{m+1}(\balg_1,..,\balg_m,\balg'),
\end{equation}
\begin{equation}\label{Baction4}
    \bi''_{l+m,n}(\balg'_1,..,\balg'_l,\balg_1,..,\balg_m,\balg'_{l+1},..,\balg'_{l+n}):=\qquad\qquad\qquad\qquad
\end{equation}
    $$\qquad\qquad\qquad\qquad=\underset{p>0}{\Sigma}\bi'_{l,p}(\balg'_1,..,\balg'_l,\be_{m+1}(\balg_1,..,\balg_m,(\balg'_{l+1},..,\balg'_{l+n}))).$$
Recall (definition \ref{BBAction}) that $\{\be_m\}$ denote both the initial action of $\2\Balg$ on $\Balg'$ and its extension to an action on
the cofree coalgebra, cogenerated by $\Balg'$.
\begin{lemma}\label{BBExtensions}
    Equations (\ref{Baction1})-(\ref{Baction4}) define the structure of a $\Bi$-algebra on $\2(\Balg'\oplus\Balg)$.
\end{lemma}
\textbf{Proof:}
    This lemma is proved in essentially the same way as lemma \ref{Extensions}. The only difference is that in the proof of
    lemma \ref{Extensions} we have to substitute $\aalg\in\Aalg$ with $\balg'_1\tensor...\tensor\balg'_m$.
$\blacksquare$

In definition of modules, algebras and coalgebras over $\Bi$-algebras we have always had actions from the left. However, actions from the right
are common, and all definitions and lemmas above, with signs suitably adjusted, apply to the case of the right action.

There will be a situation, where we have two $\Bi$-algebras $(\2\Balg,\{\bi_{m,n}\},\{\di_m\})$ and $(\2\Balg',\{\bi'_{m,n}\},\{\di'_m\})$
acting on one associative algebra $(\Aalg,\cdot,\Da)$, through $\{\be_m\}$ and $\{\be'_m\}$ respectively (definition \ref{AlgebraCoalgebra}).
The first one will act from the left and the other from the right. We will want to extend the direct product of $\2\Balg$ and $\2\Balg'$ by
$\Aalg$. It can be done if the actions of $\2\Balg$ and $\2\Balg'$ on $\Aalg$ commute, that is if for all $\balg_i\in\Balg$, $\balg'_i\in\Balg'$
and $\aalg\in\Aalg$ we have
\begin{equation}\label{Commutation}
    \be_m(\balg_1,..,\balg_{m-1},\be'_n(\aalg,\balg'_1,..,\balg'_{n-1}))=\be'_n(\be_m(\balg_1,..,\balg_{m-1},\aalg),\balg'_1,..,\balg'_{n-1}).
\end{equation}

We do it in two steps. First we extend $\2\Balg'$ by $\Aalg$ (lemma \ref{Extensions}). Then we note that, because actions of $\2\Balg$ and
$\2\Balg'$ on $\Aalg$ commute, $\{\be_m\}$ and the trivial action of $\2\Balg$ on $\2\Balg'$ induce an action of $\2\Balg$ on
$\2(\Aalg\oplus\Balg')$ (the latter considered as a $\Bi$-algebra), hence we can extend $\2\Balg$ by $\2(\Aalg\oplus\Balg')$ (lemma
\ref{BBExtensions}). The resulting $\Bi$ algebra $\2(\Aalg\oplus\Balg'\oplus\Balg)$ has the following operations
    \begin{equation}\label{2-3-1}
        \quad\bi''_{1,1}(\aalg_1,\aalg_2):=\aalg_1\cdot\aalg_2,\quad\di''_1(\aalg):=\Da(\aalg),
    \end{equation}
    \begin{equation}\label{2-3-2}
        \bi''_{m,n}(\balg_1,..,\balg_{m+n}):=\bi_{m,n}(\balg_1,..,\balg_{m+n}),\quad\di''_m(\balg_1,..,\balg_m):=\di_m(\balg_1,..,\balg_m),
    \end{equation}
    \begin{equation}\label{2-3-3}
        \bi''_{m,n}(\balg'_1,..,\balg'_{m+n}):=\bi'_{m,n}(\balg'_1,..,\balg'_{m+n}),\quad\di''_m(\balg'_1,..,\balg'_m):=\di'_m(\balg'_1,..,\balg'_m),
    \end{equation}
    \begin{equation}\label{2-3-4}
        \bi''_{m,1}(\balg_1,..,\balg_m,\aalg):=\be_{m+1}(\balg_1,..,\balg_m,\aalg),
    \end{equation}
    \begin{equation}\label{2-3-5}
        \bi''_{m,1}(\aalg_1,\balg_1,..,\balg_{m-1},\aalg_2):=\aalg_1\cdot\be_m(\balg_1,..,\balg_{m-1},\aalg_2),
    \end{equation}
    \begin{equation}\label{2-3-6}
        \bi''_{1,m}(\aalg,\balg'_1,..,\balg'_m):=\be'_{m+1}(\aalg,\balg'_1,..,\balg'_m),
    \end{equation}
    \begin{equation}\label{2-3-7}
        \bi''_{1,m}(\aalg_1,\balg'_1,..,\balg'_{m-1},\aalg_2):=\be'_m(\aalg_1,\balg'_1,..,\balg'_{m-1})\cdot\aalg_2,
    \end{equation}
    \begin{equation}\label{2-3-8}
        \bi''_{m,n}(\aalg_1,\balg_1,..,\balg_{m-1},\balg'_1,..,\balg'_{n-1},\aalg_2):=\qquad\qquad\qquad\qquad
    \end{equation}
    $$\qquad\qquad\qquad\qquad=(-1)^\sigma\be'_n(\aalg_1,\balg'_1,..,\balg'_{n-1})\cdot\be_m(\balg_1,..,\balg_{m-1},\aalg_2),$$
where $\balg_i\in\Balg$, $\balg'_i\in\Balg'$, $\aalg_i\in\Aalg$, $\sigma=(\balg_1+...+\balg_{m-1})(\balg'_1+...+\balg'_{n-1})$ and for the rest
of possibilities the values are $0$.

\section{Deformations of morphisms}\label{SectionMorphisms}
\subsection{$\G$-structure}\label{SubsectionGstructure}$\quad$\\
In this subsection we define $\Bi$-structure on the deformation complex of a morphism of associative algebras. Then existence of a
$\G$-structures follows from the well known fact (\cite{TT} section 3, \cite{Hin2} sections 6,7) that all $\Bi$-algebras are also $\G$-algebras
(in a non-unique way).

The $\Bi$-structure is constructed as an extension of the direct product of the Hochschild complexes of the two associative algebras (subsection
\ref{Ext}). As usual we start with just a pair of graded $\K$-spaces $\Aone$, $\Atwo$, and produce a $\Bi$-algebra. Then using results of the
previous section we deform this $\Bi$-structure by a solution of the Maurer-Cartan equation, that corresponds to the associative structures on
$\Aone$, $\Atwo$ and the morphism between them, that we wish to deform.

Let $\Aone$, $\Atwo$ be graded $\K$-vector spaces. Define
    $$\GOne:=\SM\Hom((\s\Aone)\tensorM,\s\Aone),\quad\GTwo:=\SM\Hom((\s\Atwo)\tensorM,\s\Atwo),$$
    $$\GOneTwo:=\SM\Hom((\s\Aone)\tensorM,\s\Atwo).$$
We have $\Bi$-structures on $\2\GOne$ and $\2\GTwo$ (\cite{Get} subsection 5.2). The operations are as follows
    $$\bi_{1,m}(\gTwo_1,..,\gTwo_{m+1}):=\gTwo_1\circ(\gTwo_2\tensor..\tensor\gTwo_{m+1}),$$
    $$\bi_{1,m}(\gOne_1,..,\gOne_{m+1}):=\gOne_1\circ(\gOne_2\tensor..\tensor\gOne_{m+1}),$$
where $\gOne_i\in\GOne$, $\gTwo_i\in\GTwo$, $\gOne_1\in\Hom((\s\Aone)^{\tensor^{\geq m}},\s\Aone)$, $\gTwo_1\in\Hom((\s\Atwo)^{\tensor^{\geq
m}},\s\Atwo)$, and $\{\di_m\}$ operations are trivial. We are going to extend these $\Bi$-algebras.

First we consider the cofree coalgebra, cogenerated by $\GOneTwo$, we will denote it by $\TG$. We have two actions on $\TG$: $\2\GOne$ acts from
the right and $\2\GTwo$ from the left (definition \ref{AlgebraCoalgebra}). These actions are defined by the corresponding $\Bi$-morphisms (lemma
\ref{Representation}) as follows:
    $$\CacR_2(\gOneTwo_1\tensor..\tensor\gOneTwo_n,\gOne):=\underset{i=1}{\overset{n}{\Sigma}}(-1)^{\sigma'(i)}\gOneTwo_1\tensor..\tensor
    \gOneTwo_i\circ\gOneHat\tensor..\tensor\gOneTwo_n,$$
    $$\CacL_2(\gTwo,\gOneTwo_1\tensor..\tensor\gOneTwo_n):=\underset{i=0}{\overset{n-m}{\Sigma}}(-1)^{\sigma(i)}\gOneTwo_1\tensor..\tensor
    \gTwo\circ(\gOneTwo_{i+1}\tensor..\tensor\gOneTwo_{i+m})\tensor..\tensor\gOneTwo_n,$$
where $\gOneTwo_i\in\GOneTwo$, $\gOne\in\GOne$, $\GTwo\in\Hom(\Atwo\tensorM,\Atwo)$, $\sigma'(i)=\gOne(\gOneTwo_{i+1}+...+\gOneTwo_n)$,
$\sigma(i)=\gTwo(\gOneTwo_1+...+\gOneTwo_i)$, and $\gOneHat$ stands for the coderivation on $\TC(\Aone)$, generated by $\gOne$. Direct
calculation shows that $\CacR$ and $\CacL$ indeed define actions of $\Bi$-algebras $\2\GOne$, $\2\GTwo$ on $\TG$.

Since compositions from the right and from the left obviously commute, we see that the actions of $\2\GOne$ and $\2\GTwo$ on $\TG$ commute, i.e.
equations (\ref{Commutation}) are satisfied. Applying lemma \ref{Cobar} we have actions of $\2\GOne$ and $\2\GTwo$ on the cobar construction of
$\TG$. We will denote this algebra by $(\GG,\de)$, where $\GG=\Om(\2\TG)\tensorM$ and $\de$ is given by its values on the generators of $\GG$ as
follows:
\begin{equation}\label{CobarEquation}
    \de(\2(\gOneTwo_1,..,\gOneTwo_n))=\underset{i=1}{\overset{n-1}{\Sigma}}(-1)^{\sigma(i)}\2(\gOneTwo_1,..,\gOneTwo_i)\tensor
    \2(\gOneTwo_{i+1},..,\gOneTwo_n),
\end{equation}
where $\sigma(i)=\gOneTwo_1+...+\gOneTwo_i+1$, $\gOneTwo_i\in\GOneTwo$.

Clearly the actions of $\2\GOne$ and $\2\GTwo$ on $(\GG,\de)$ also
commute. Therefore we can define a $\Bi$-structure on
$\2(\GG\oplus\GOne\oplus\GTwo)$. However, this $\Bi$-algebra would
not yet be the correct one to describe deformations of morphisms.
We have to modify it in two respects.

The first modification is that we have to change sign of the action of $\2\GOne$ on $\GG$. The reason for this is that the defining equation of
an $\A$-morphism between $\A$-algebras is
    $$\widetilde{\gamma}\circ\widehat{\alpha}-\widehat{\beta}\circ\widetilde{\gamma}=0,$$
where $\alpha$, $\beta$ are elements of $\GOne$, $\GTwo$ respectively, $\gamma\in\GOneTwo$, hat and tilde denote respectively the coderivations
and the coalgebra morphisms, generated by the cochains. That minus sign before the second summand on the l.h.s. is what requires the
modification of the action. So we change the action of $\2\GOne$ on $\GG$ as follows:
    $$\CacR_m(\2\gOneTwo,\gOne_1,..,\gOne_{m-1}):=(-1)^{\gOne_1+..+\gOne_{m-1}}\2(\gOneTwo\circ(\gOne_1\tensor..\tensor\gOne_{m-1})).$$
Clearly actions of $\2\GOne$ and $\2\GTwo$ on $\GG$ still commute, hence we have a $\Bi$-structure on $\2(\GG\oplus\GOne\oplus\GTwo)$, that we
will denote by $\{\bi_{m,n}\}$, $\{\di_m\}$. These operations can be read from equations (\ref{2-3-1})-(\ref{2-3-8}). Explicitly
\begin{equation}\label{FinalStructure1}
    \bi_{1,1}(\2(\gOneTwo_1,..,\gOneTwo_m),\2(\gOneTwo'_1,..,\gOneTwo'_n)):=\2(\gOneTwo_1,..,\gOneTwo_m)\tensor\2(\gOneTwo'_1,..,\gOneTwo'_n),
\end{equation}
\begin{equation}
    \bi_{1,m}(\gOne_1,..,\gOne_{m+1})=\gOne_1\circ(\gOne_1\tensor...\tensor\gOne_{m+1}),
\end{equation}
\begin{equation}
    \bi_{1,m}(\gTwo_1,..,\gTwo_{m+1}):=\gTwo_1\circ(\gTwo_2\tensor..\tensor\gTwo_{m+1}),
\end{equation}
\begin{equation}
    \bi_{m,1}(\gTwo_1,..,\gTwo_m,\2(\gOneTwo_1,..,\gOneTwo_n))=\2((\gTwo_1\tensor..\tensor\gTwo_m)\circ(\gOneTwo_1\tensor..\tensor\gOneTwo_n)),
\end{equation}
\begin{equation}
    \bi_{1,m}(\2(\gOneTwo_1,..,\gOneTwo_n),\gOne_1,..,\gOne_m)=
\end{equation}
    $$(-1)^{\gOne_1+..\gOne_m}\2((\gOneTwo_1\tensor..\tensor\gOneTwo_n)\circ(\gOne_1\tensor..\tensor\gOne_m)),$$
\begin{equation}
    \bi_{m+1,1}(\2(\gOneTwo_1,..,\gOneTwo_l),\gTwo_1,..,\gTwo_m,\2(\gOneTwo'_1,..,\gOneTwo'_n))=
\end{equation}
    $$=\2(\gOneTwo_1,..,\gOneTwo_l)\tensor\2((\gTwo_1\tensor..\tensor\gTwo_m)\circ(\gOneTwo'_1\tensor..\tensor\gOneTwo'_n)),$$
\begin{equation}
    \bi_{1,m+1}(\2(\gOneTwo_1,..,\gOneTwo_n),\gOne_1,..,\gOne_m,\2(\gOneTwo'_1,..,\gOneTwo'_l))=
\end{equation}
    $$=(-1)^{\gOne_1+..+\gOne_m}\2((\gOneTwo_1\tensor..\tensor\gOneTwo_n)\circ(\gOne_1\tensor..\tensor\gOne_m))\tensor
    \2(\gOneTwo'_1,..,\gOneTwo'_l),$$
\begin{equation}\label{FinalStructure2}
    \bi_{m+1,n+1}(\2(\gOneTwo_1,..,\gOneTwo_p),\gTwo_1,..,\gTwo_m,\gOne_1,..,\gOne_n,\2(\gOneTwo'_1,..,\gOneTwo'_q))=
\end{equation}
    $$=(-1)^\sigma\2((\gOneTwo_1,..,\gOneTwo_p)\circ(\gOne_1\tensor..\tensor\gOne_n))\tensor\2((\gTwo_1\tensor..\tensor\gTwo_m)\circ
    (\gOneTwo'_1\tensor..\tensor\gOneTwo'_q)),$$
where $\gOne\in\GOne$, $\gTwo\in\GTwo$, $\gOneTwo\in\GOneTwo$, $\sigma=(\gOne_1+...+\gOne_n)(\gTwo_1+...+\gTwo_m+1)$. For the rest of
$\2(\GG\oplus\GOne\oplus\GTwo)$, i.e. for elements of tensor powers of $\2\TG$, the operations are easily derived from these equations, since
$\GG$ is an almost free algebra.

The second modification is that we have to complete $\GG$ with respect to a certain filtration. The reason for it, is that our representation of
morphisms will consist of infinite series of multi-linear maps, rather than finite sums of them.

These infinite series are of the form $\Sm(\gOneTwo)\tensorM$ ($\gOneTwo\in\GOneTwo$), i.e. we have to complete $\TG$ with respect to the
grading by tensor powers of $\GOneTwo$. In turn, differential on the cobar construction $\GG$ requires then completion of $\2\TG\tensor\2\TG$,
and so on.

Denote $\GGP:=(\2\TG)^{\tensor^p}$, and by $\BGGP$ the corresponding completion. Then $\BGG:=\underset{p>0}{\bigoplus}\BGGP$ is a differential
associative algebra, since multiplication and differential on $\GG$ are continuous with respect to the grading above. Similarly actions of
$\2\GOne$ and $\2\GTwo$ on $\GG$ extend to actions on $\BGG$, and we have a $\Bi$-structure on $\2(\BGG\oplus\GOne\oplus\GTwo)$.

Now we have to show that $\2(\BGG\oplus\GOne\oplus\GTwo)$ is the correct $\Bi$-algebra to describe deformations of associative structures on
$\Aone$, $\Atwo$ and a morphism between them, i.e. we have to analyze the underlying dg Lie algebra.

\subsection{Lie structure and deformations}\label{SubsectionLieStructure}$\quad$\\
In this section we prove that $\2(\BGG\oplus\GOne\oplus\GTwo)$ is the correct $\Bi$-algebra to describe deformations of morphisms between
non-positively graded dg associative algebras.

We do that by first showing that solutions of the Maurer-Cartan equation in the underlying dg Lie algebra are in bijective correspondence with
$\A$-structures on $\Aone$, $\Atwo$ (i.e. with structures of $\A$-algebras on $\Aone$, $\Atwo$ and an $\A$-morphism between them).

Then we use lemma \ref{BDeformation} and show that the deformations of the $\Bi$-structure on $\2(\BGG\oplus\GOne\oplus\GTwo)$, given by
solutions of the Maurer-Cartan equation, have the correct underlying dg Lie algebras to describe deformations of cobar constructions of the
corresponding $\A$-structures.

To do that we prove that if we have two pairs of graded
$\K$-spaces $\Aone$, $\Atwo$ and $\Aone'$, $\Atwo'$, and two
solutions of the Maurer-Cartan equations in the corresponding
$\Bi$-algebras $\2(\BGG\oplus\GOne\oplus\GTwo)$ and
$\2(\BGG'\oplus\GOne'\oplus\GTwo')$, s.t. cobar constructions of
the $\A$-structures, given by the solutions, are quasi-isomorphic,
the deformations of $\2(\BGG\oplus\GOne\oplus\GTwo)$ and
$\2(\BGG'\oplus\GOne'\oplus\GTwo')$, given by these solutions, are
quasi-isomorphic $\Bi$-algebras. This allows us to assume that the
morphism of associative algebras we wish to deform is injective,
which makes the dg Lie algebra easy to handle.

\subsubsection{Structures over a point}$\quad$\\
It is well known that given a $\Bi$-algebra $(\2\Balg,\{\bi_{m,n}\},\{\di_m\})$, we have a dg Lie algebra $(\Balg,[,],\di_1)$, with the bracket
defined as follows:
\begin{equation}\label{CobarDifferential}
    [\balg_1,\balg_2]=\bi_{1,1}(\balg_1,\balg_2)+(-1)^{\balg_1\balg_2+1}\bi_{1,1}(\balg_2,\balg_1).
\end{equation}
Applying this to the $\Bi$-structure on $\2(\BGG\oplus\GOne\oplus\GTwo)$ we get a dg Lie algebra, that we will denote by $\GMorphism$. As
the following proposition shows, $\GMorphism$ is the correct dg Lie algebra to describe deformations over a point, even if $\Aone$, $\Atwo$ are
not non-positively graded.
\begin{proposition}\label{MorPoint}
    There is a bijection between the set of $\A$-structures on $\Aone,\Atwo$ (i.e. $\A$-algebras on $\Aone,\Atwo$, and an $\A$-morphism
    $\Aone\rightarrow\Atwo$) and the set of solutions of the Maurer-Cartan equation in $\GMorphism$.
\end{proposition}
\textbf{Proof:}
    Let $\gOne$, $\gTwo$, $\gOneTwo$ be an $\A$-structure on $\Aone$, $\Atwo$, i.e. $\gOne$, $\gTwo$ are codifferentials of degree 1 on
    $\B(\Aone)$, $\B(\Atwo)$ respectively and $\gOneTwo$ is a degree 0 coalgebra morphism $\B(\Aone)\rightarrow\B(\Atwo)$, where $\B$ stands
    for the bar construction functor. Here we use the fact that
    $\B(\Aone)$, $\B(\Atwo)$ are almost cofree coalgebras and therefore codifferentials on them and morphisms between them can be represented by their
    corestrictions to cogenerators, i.e. $\gOne\in\GOne$, $\gTwo\in\GTwo$, $\gOneTwo\in\GOneTwo$.

    The condition on $\gOne$, $\gTwo$, $\gOneTwo$ to be an $\A$-structure is
    \begin{equation}\label{MorphismConditions}
        \gOneHat^2=0,\quad\gTwoHat^2=0,\quad\gOneTwo\circ\gOneHat=\gTwo\circ\gOneTwoTilde,
    \end{equation}
    where $\gOneHat$ is the coderivation cogenerated by $\gOne$, and $\gOneTwoTilde$ is the coalgebra morphism, cogenerated by
    $\gOneTwo$.

    Define an element $\GMe\in\GMorphism$ by $\GMe:=\gOne+\gTwo+\Sm\2(\gOneTwo)\tensorM$. By definition of the cobar construction, $\GG$ is an
    almost free algebra, generated by $\2\TG$, and has a differential given by equation (\ref{CobarEquation}). Because degree of
    $\gOneTwo$ is $0$, from equation (\ref{CobarDifferential}) and equations (\ref{FinalStructure1})-(\ref{FinalStructure2}) we have
        $$\de(\Sm\2(\gOneTwo)\tensorM)+\frac{1}{2}[\Sm\2(\gOneTwo)\tensorM,\Sm\2(\gOneTwo)\tensorM]=0.$$
    The rest of the Maurer-Cartan equation for $\GMe$ is
        $$[\gOne,\gOne]=0,\quad[\gTwo,\gTwo]=0,\quad[\gOne+\gTwo,\Sm\2(\gOneTwo)\tensorM]=0,$$
    but these are exactly equations (\ref{MorphismConditions}), so $\GMe$ is a solution.

    Now let $\GMe$ be a solution of the Maurer-Cartan equation in $\GMorphism$. Since $\GMorphism=\BGG\oplus\GOne\oplus\GTwo$, we can
    write
        $$\GMe=\gOne+\gTwo+\underset{1\leq p\leq q}{\Sigma}\gap,$$
    where $\gOne\in\GOne$, $\gTwo\in\GTwo$, $\gap\in\BGGP$. Suppose $q>1$, then $[\gaq,\gaq]\in\BGGQQ$, and from the
    Maurer-Cartan equation we conclude $[\gaq,\gaq]=0$, but since this bracket comes from commutator in a free associative algebra,
    and $\Deg(\gaq)=1$, we have $\gaq=0$, i.e. $q=1$. Then again from the Maurer-Cartan equation it follows that
        $$\GMe=\gOne+\gTwo+\Sm\2(\gOneTwo)\tensorM,$$
    for some $\gOneTwo\in\GOneTwo$, and Maurer-Cartan equation for $\GMe$ translates into equations (\ref{MorphismConditions}) for
    $\gOne$, $\gTwo$, $\gOneTwo$, i.e. they constitute an $\A$-structure.
$\blacksquare$

\subsubsection{Invariance with respect to quasi-isomorphisms}$\quad$\\
As usual, having a dg Lie algebra that controls all possible structures over a point, to describe deformations of a given structure one has to
take the corresponding solution of the Maurer-Cartan equation and deform with it the differential in the dg Lie algebra. This is done in the
following definition.
\begin{definition}\label{DefinitionLie}
    Let $\go:=\goo+\got+\gaz$ be a solution of the Maurer-Cartan equation in $\GMorphism$, where $\goo\in\GOne$,
    $\got\in\GTwo$, $\gaz=\Sm\2(\goot)\tensorM$, $\goot\in\GOneTwo$. Let $\mor$ be the corresponding
    $\A$-morphism of $\A$-algebras (proposition \ref{MorPoint}). We define \underline{dg Lie algebra $\GDefMor$} to be the same graded Lie algebra as
    $\GMorphism$, but with the differential being
        $$\Dg:=\de+[\go,-],$$
    where $\de$ is given by equation (\ref{CobarEquation}). Also we define \underline{$\Bi$-algebra $\2\BDefMor$} to be the deformation of
    $\2(\BGG\oplus\GOne\oplus\GTwo)$ by $\go$ (lemma \ref{BDeformation}).
\end{definition}
It is easy to see that $\GDefMor$ is the underlying dg Lie algebra of $\2\BDefMor$. Therefore, if we prove that $\GDefMor$ is the correct dg Lie
algebra to describe deformations of $\mor$, we would show that $\2\BDefMor$ is the correct $\Bi$-algebra.

To prove this, first we have to show that by choosing another morphism $\morr$, that is quasi-isomorphic to $\mor$, we get a $\Bi$-algebra
$\2\BDefMorr$, that is quasi-isomorphic to $\2\BDefMor$, where by quasi-isomorphic $\Bi$-algebras we mean algebras, s.t. there is a chain of
quasi-isomorphisms connecting them, with quasi-isomorphism being the usual notion for operadic algebras. This will clearly imply, that the
corresponding dg Lie algebras are quasi-isomorphic as well.

Then we can choose a representative of the quasi-isomorphism class of $\mor$, that is suitable for the proof. Recall that a morphism between
morphisms is a pair of morphisms, that make up a commutative square. A quasi-isomorphism between morphisms is a pair of quasi-isomorphisms.
\begin{lemma}\label{AlgebraChange}
    If $\mor$ and $\morr$ are quasi-isomorphic, $\2\BDefMor$ and $\2\BDefMorr$ are quasi-isomorphic as well.
\end{lemma}
\textbf{Proof:}
    Let $\mor:\Aone\rightarrow\Atwo$ and $\morr:\Aone\rightarrow\Btwo$ be two morphisms of associative algebras. Let
    $\phi:\Atwo\rightarrow\Btwo$ be a quasi-morphism, such that $\morr=\psi\circ\mor$. Define
        $$\HTwo:=\SM\Hom((\s\Btwo)\tensorM,\s\Btwo),\quad\HOneTwo:=\SM\Hom((\s\Aone)\tensorM,\s\Btwo),$$
        $$\Phi:=\SM\Hom((\s\Atwo)\tensorM,\s\Btwo).$$
    Let $\BHH$ be the completion of the cobar construction of $\THH$ (as in definition of $\BGG$).

    We have two $\Bi$-algebras: $\2(\BGG\oplus\GOne\oplus\GTwo)$ and $\2(\BHH\oplus\GOne\oplus\HTwo)$. We deform them, according to lemma
    \ref{BDeformation}, by solutions of the Maurer-Cartan equation, that correspond to $\mor$ and $\morr$. The results are two $\Bi$-algebras
    $\2\BDefMor$ and $\2\BDefMorr$. Now we construct a chain of quasi-isomorphism, connecting these two $\Bi$-algebras.

    Since $\phi$ is a quasi-isomorphism,
    composition with it defines a quasi-isomorphism $\phi_*:\BGG\rightarrow\BHH$, where differentials on $\BGG$ and $\BHH$ are
    restrictions from $\BDefMor$ and $\BDefMorr$. It also defines two maps from $\GTwo$ and $\HTwo$
    to $\Phi$. We will denote by $\Ha$ the vector space that is the fiber product of these two maps.

    We have $\Bi$ structures on $\2\GTwo$ and $\2\HTwo$. Elements of $\Ha$ are pairs of elements of
    $\GTwo$ and $\HTwo$, that satisfy a relation. It is easily seen, that componentwise application of the $\Bi$-operations to such pairs
    satisfies this relation too. Therefore $\2\Ha$ is a $\Bi$-algebra. Since it is a fiber product, it is mapped into $\2\GTwo$ and $\2\HTwo$,
    and since $\phi$ is a quasi-isomorphism, these maps are quasi-isomorphisms as well.

    Actions of $\2\GTwo$, $\2\HTwo$ on $\BGG$, $\BHH$ respectively induce actions of $\2\Ha$ on them, and $\phi_*$ induces a morphism between
    these actions. It is a quasi-isomorphism, because $\phi_*$ is one.

    In total we have a string of quasi-isomorphisms of $\Bi$-algebras
        $$\2\BDefMor\leftarrow\2(\BGG\oplus\GOne\oplus\Ha)\rightarrow\2(\BHH\oplus\GOne\oplus\Ha)\rightarrow\2\BDefMorr.$$
    The case of changing the domain of $\mor$ is done similarly.
$\blacksquare$

\subsubsection{Deformations over dg Artin algebras}$\quad$\\
Deformations of $\mor$ as an $\A$-morphism are mapped (by means of the cobar construction) to deformations of the cobar-bar construction of
$\mor$. These are described by the simplicial groupoid $\DEFf$ (\cite{Hin1},\cite{Bor}(definition 3)).

Let $\dgart$ be the category of local dg Artin algebras over $\K$, let $\R\in\dgart$. Objects of $\DEFf(\R)$ are cofibrations between cofibrant
dg associative algebras over $\R$, such that their reductions modulo the maximal ideal $\RI$ of $\R$ is $\Omega\B(\mor)$ (recall that $\Omega$
stands for the cobar construction, and $\B$ denotes the bar construction). Here we use lemma \ref{AlgebraChange}, and assume that $\mor$ is
injective, which implies that its cobar-bar construction is a cofibration.

The classical deformation groupoid $\Deff(\R)$ consists of the same objects as $\DEFf(\R)$, but its sets of morphisms are sets of connected
components of the mapping spaces in $\DEFf(\R)$.

The question whether $\GDefMor$ (definition \ref{DefinitionLie}) is the correct dg Lie algebra amounts to whether the simplicial Deligne
groupoid $\DELg$ (\cite{Hin1} subsection 3.1) corresponding to $\GDefMor$ is weakly equivalent to $\DEFf$. For morphisms between non-positively
graded algebras this is proved in the following theorem.
\begin{theorem}\label{NonPositively}
    Let $\mor$ be a morphism between non-positively graded dg associative algebras. There is a weak equivalence of simplicial
    categories
        $$\FR:\DELg(\R)\rightarrow\DEFf(\R),$$
    that is natural in $\R$.
\end{theorem}
\textbf{Proof:}
    We proceed as follows. First we construct a dg Lie subalgebra $\Ha$ of $\GDefMor$, whose simplicial Deligne groupoid is easier to connect with the
    deformation groupoid $\DEFf(\R)$.

    Then we show that inclusion of this subalgebra in $\GDefMor$ is a quasi-isomorphism, and hence these two dg Lie
    algebras have weakly equivalent simplicial Deligne groupoids (\cite{Hin1}, corollary 3.3.2).

    Next we construct a functor $\Fr$ from the classical Deligne groupoid $\Del(\Ha)(\R)$, corresponding to the subalgebra, to the classical deformation
    groupoid $\Deff(\R)$. We show that this functor is an equivalence of categories.

    Finally we use the fact that nerves of $\DELg(\R)$ and $\DEFf(\R)$ are bisimplicial sets,
    consisting of nerves of certain sub-groupoids of the classical Deligne groupoid $\Delg(\R)$ and $\Deff(\R)$ respectively. This implies weak
    equivalence of the bisimplicial sets and we are done.

    \textbf{1. Construction of the dg Lie subalgebra:} Let $\go=\goo+\got+\gaz$, where $\goo\in\GOne$, $\got\in\GTwo$, $\gaz\in\BGG$, be the solution
    corresponding to $\mor$ according to proposition
    \ref{MorPoint}. Let $\Ha$ be the $\K$-subspace of $\GDefMor$, defined as follows
        $$\Ha:=\{\Sh\in\GOne\oplus\GTwo\text{ s.t. }\Dg(\Sh)\in\GOne+\GTwo\}.$$

    We claim that $\Ha$ is a dg Lie subalgebra of $\GDefMor$. Indeed, $\Sh\in\Ha$ if and only if $[\gaz,\Sh]=0$, therefore, since in that case
        $$\Dg(\Sh)=[\goo+\got,\Sh],$$
    and in general $[\gaz,\goo+\got]=0$, we have that $[\gaz,\Dg(\Sh)]\in\GOne\oplus\GTwo$, i.e. $\Ha$ is
    a subcomplex of $\GDefMor$. If $\Sh$, $\Sh'\in\Ha$, then, since $\GOne+\GTwo$ is a graded Lie subalgebra,
        $$\Dg([\Sh,\Sh'])=[\Dg(\Sh),\Sh']+(-1)^\Sh[\Sh,\Dg(\Sh')]$$
    also belongs to $\GOne+\GTwo$, i.e. $\Ha$ is closed under the bracket.

    \textbf{2. The inclusion $\Ha\hookrightarrow\GDefMor$ is a quasi-isomorphism:} First we prove it is surjective on cohomology.
    Let $\GMe\in\GDefMor$ be a cocycle for $\Dg$. We write $\GMe=\gOne+\gTwo+\underset{1\leq p\leq q}{\Sigma}\gap$, where $\gOne\in\GOne$,
    $\gTwo\in\GTwo$, $\gap\in\BGGP$. Since $\Dg(\GMe)=0$, we have $\de(\gaq)+[\gaz,\gaq]=0$. If $\gaq$ was also a coboundary for $\de+[\gaz,-]$,
    then $\GMe$ would have been cohomologous to an element of $\GOne\oplus\GTwo\oplus(\underset{1\leq p\leq q-1}{\bigoplus}\BGGP)$. Hence the following
    lemma.
    \begin{lemma}\label{Surprize}
        Every cocycle of $\de+[\gaz,-]$ in $\underset{p\geq 2}{\bigoplus}\BGGP$ is also a coboundary.
    \end{lemma}
    \textbf{Proof:}
        Let $K$ be the kernel of the projection $\BGG\rightarrow\2\GOneTwo$. It is graded by $K_p:=K\cap\BGGP$, and filtered by tensor powers
        of $\GOneTwo$. From definitions of $\de$ and $[\gaz,-]$ it follows that
        $(K,\de+[\gaz,-])$ is a filtered complex. The filtration is obviously exhaustive. It is also weakly convergent because if we denote by
        $F^nK$ the $n$-th part of the filtration, then $F^nK/F^{n+r}K$ is acyclic for all $n,r\in\mathbb N$ (to see it note that the filtration on
        $(F^nK/F^{n+r}K,\de+[\gaz,-])$ is bounded and $H^*(F^nK/P^{n+r}K,\de)=0$). Therefore the associated spectral
        sequence converges, but its first term is trivial, so $K$ is acyclic. Hence every cocycle in $\BGG$ is cohomologous to one in $\BGGone$,
        but there are no
        coboundaries with non-trivial projections on $\BGGone$, therefore any cocycle in $\underset{p\geq 2}{\bigoplus}\BGGP$ is a coboundary.
    $\blacksquare$

    From lemma \ref{Surprize} it follows that any cocycle in $\GDefMor$ is cohomologous to an element of $\GOne\oplus\GTwo\oplus\BGGone$.
    Let $\GMe$ be such. Then we can write $\GMe=\gOne+\gTwo+\Sm\2\gam$, where
    $\gOne\in\GOne$, $\gTwo\in\GTwo$, $\gam\in\GOneTwo\tensorM$. Since $\GMe$ is a cocycle we have
        $$\de(\2\gat)+\2\goot\tensor\2\gao+(-1)^\gao\2\gao\tensor\2\goot=0,$$
    where $\gaz=\Sm\2(\goot)\tensorM$ (recall that $\Deg(\goot)=0$). Therefore
        $$\gat=\goot\tensor\gao+\gao\tensor\goot.$$
    Applying $\Dg$ again we find that
        $$\gath=\goot\tensor\goot\tensor\gao+\goot\tensor\gao\tensor\goot+\gao\tensor\goot\tensor\goot$$
    and so on. We see that all $\Zh_m$ are determined by $\gao$, in particular if $\gao=0$, then $\GMe\in\GOne\oplus\GTwo$.

    Now we use lemma \ref{AlgebraChange}, from which it follows we may suppose, that $\mor$ is an injective morphism of associative algebras.
    Consider a
    decomposition $\Atwo=U\oplus V$, where $U$ is the image of $\mor$, and define $\mor^{-1}$ to be the inverse of $\mor$ on $U$ and 0 on $V$. Let
    $\psi\in\Hom((\s\Aone)\tensorM,\s\Atwo)$, define $\beta\in\Hom((\s\Atwo)\tensorM,\s\Atwo)$ by
        $$\beta(b_1...b_m):=\psi(\mor^{-1}(b_1)...\mor^{-1}(b_m)).$$
    If we do this with $\gao$ we have $\GMe-\Dg(\beta)\in\GOne\oplus\GTwo$. So every cocycle in $\GDefMor$ is cohomologous to one in
    $\GOne\oplus\GTwo$ and those are obviously elements in $\Ha$.

    Therefore the inclusion $\Ha\hookrightarrow\GDefMor$ is surjective on cohomology. To see that it is injective on cohomology
    note that if $\GMe\in\GDefMor$ and $\Dg(\GMe)\in\Ha$, then obviously $\Dg(\GMe)\in\GOne\oplus\GTwo$. It may happen that the projection of $\GMe$ on
    $\BGG$ is not 0. In that case proceeding as in the proof of lemma \ref{Surprize} we can find an $\GMe'\in\GOne\oplus\GTwo\oplus\BGGone$, s.t.
    $\Dg(\GMe')=\Dg(\GMe)$. Then projection of $\GMe'$ on $\BGGone$ is a cocycle for $\de+[\gaz,-]$, and using injectivity of $\mor$ as above we can find
    $\GMe''\in\GOne\oplus\GTwo$, s.t. $\Dg(\GMe'')=\Dg(\GMe)$. Then obviously $\GMe''\in\Ha$. So a cocycle for $\Dg$ in $\Ha$ is a coboundary if
    and only if it is such in $\GDefMor$, i.e. the inclusion is injective on cohomology.

    We have constructed a dg Lie algebra $\Ha$ that is quasi-isomorphic to $\GDefMor$. According to \cite{Hin1} corollary 3.3.2 the corresponding
    simplicial
    Deligne groupoids $\DEL(\Ha)$ and $\DEL(\GDefMor)$ are weakly equivalent. Now we show that $\DEL(\Ha)$ is weakly equivalent to $\DEFf$.

    \textbf{3. Construction of $\Fr$:} Let $\R\in\dgart$. Let $\GMe\in\Ha\tensor\RI$ be a solution of the Maurer-Cartan equation in $\Ha\tensor\RI$.
    We can write
    $\GMe=\gOne+\gTwo$, where $\gOne\in\GOne\tensor\RI$ and $\gTwo\in\GTwo\tensor\RI$. From the dg Lie structure on $\Ha$ it is clear that
    $\gOne$, $\gTwo$ determine deformations of $\A$-algebras $\Aone$, $\Atwo$ respectively. In addition they satisfy
        $$[\gaz,\gOne]+[\gaz,\gTwo]=0.$$
    From definition of the bracket one sees that this equation is equivalent to following one
        $$\gTwoHat\circ\mor=\mor\circ\gOneHat,$$
    where as before, $\gOneHat$, $\gTwoHat$ are coderivations on $\B(\Aone)$, $\B(\Atwo)$, cogenerated by $\gOne$, $\gTwo$ respectively, and we
    use the same symbol $\mor$ for the $\R$-linear extension of $\mor$. So every solution in $\g\tensor\RI$ represents a
    deformation of the $\A$-structure $\mor$.

    Equivalences between solutions in $\Ha\tensor\RI$ are exponentials of elements of degree 0 in $\Ha\tensor\RI$. These exponentials represent
    morphisms between the corresponding $\A$-deformations, s.t. their reductions modulo $\RI$ is the identity. So by means of the cobar
    construction we have a functor $\Fr:\Del(\Ha)(\R)\rightarrow\Deff(\R)$.

    \textbf{4. Equivalence of classical groupoids:} We claim that $\Fr$ is essentially surjective, i.e. every object in
    $\Deff(\R)$ is isomorphic to an image of $\Fr$. Indeed,
    every deformation of an associative algebra is quasi-isomorphic to cobar construction of its deformation as an $\A$-algebra.
    Similarly every object in $\Deff(\R)$ is isomorphic to cobar construction of an $\A$-deformation of $\mor$.

    Now let $(\gOne,\gTwo,\gOneTwo)$ be an $\A$-deformation of $\mor$, i.e. $\gOne\in\GOne$, $\gTwo\in\GTwo$, $\gOneTwo\in\GOneTwo$.
    Since $\mor$ is injective, proceeding as before (in the proof that $\Ha\hookrightarrow\GDefMor$ is a quasi-isomorphism), we can find an element
    $\be\in\GTwo\tensor\RI$, s.t.
    $\mor+\gOneTwoTilde=\e(\be)\circ\mor$, i.e. the deformation $(\gOne,\gTwo,\gOneTwo)$ is equivalent to
    $(\gOne,\e(\be)(\gTwo),0)$, which is in the image of $\Fr$. So $\Fr$ is essentially surjective.

    Next we show that $\Fr$ is injective on the skeleton, i.e. any two objects in $\Del(\Ha)(\R)$ are equivalent if and only if their images
    under $\Fr$ are equivalent. Let $\Fr(\Zh)$, $\Fr(\Zh')$ be objects in the image of $\Fr$. If they are equivalent in $\Deff(\R)$, their equivalence
    can be represented by $\A$-morphisms, i.e.  there are two
    $\A$-morphisms $\mu:\B(\Aone)\rightarrow\B(\Aone)$ and $\nu:\B(\Atwo)\rightarrow\B(\Atwo)$, whose cobar constructions constitute an
    equivalence from $\Fr(\Zh)$ to $\Fr(\Zh')$. Let $\al$, $\be$ be elements of $\GOne$, $\GTwo$ respectively, such that
        $$\mu=\e(\al),\quad\nu=\e(\be).$$
    Then $\Zh'=\e(\al+\be)(\Zh)$, and $\Zh$, $\Zh'$ are equivalent in $\Del(\Ha)(\R)$.

    \textbf{5. Equivalence of simplicial groupoids: } The objects in $\DEFf(\R)$ are deformations of $\mor$ over $\R$, and for any two of them
    $X$, $Y$ the mapping space $\HOM(X,Y)$ has
    the following components: $\HOM_n(X,Y)$ is the set of morphisms from the $\Omega_n$-linear extension of $X$ to that of $Y$, s.t. their
    reductions modulo $\RI$ is the identity. Here $\Omega_n$ stands for the algebra of polynomial forms on an $n$-simplex (\cite{Bou} chapter 1).
    Similarly for the mapping spaces in $\DEL(\Ha)(\R)$ (\cite{Hin1} subsection 3.1).

    Clearly $\Omega_n\tensor\RI$ is an
    Artinian algebra, and therefore $\Deff(\Omega_n\tensor\RI)$ is equivalent to $\Del(\Ha)(\Omega_n\tensor\RI)$, as we have shown above.
    In particular the
    subcategories of $\Omega_n$-linear extensions are equivalent too. Therefore their nerves are weakly equivalent simplicial sets.

    These
    simplicial sets are the components of the bisimplicial sets, that are the componentwise nerves of $\DEL(\Ha)(\R)$ and $\DEFf(\R)$. Therefore $\FR$
    induces a weak equivalence of the nerves of $\DEL(\Ha)(\R)$ and $\DEFf(\R)$ (\cite{Hir} theorems 15.11.6, 15.11.11), and hence a weak
    equivalence of the simplicial groupoids themselves (\cite{Hin1} proposition 6.3.3).
$\blacksquare$

\subsection{Diagram algebra and cohomology}\label{SubsectionDiagram}$\quad$\\
In \cite{GS1} section 23, it is shown that cohomology of a diagram of associative algebras is isomorphic to cohomology of a certain algebra,
that is built from the diagram. This algebra is called there the \underline{diagram algebra}, and we will also use this name here.

In particular this result is true for the case of a morphism of dg associative algebras $\mor:\Aone\rightarrow\Atwo$. In this case the diagram
algebra $\AD$ is as follows (\cite{GS2} section 2): as a dg $\K$-space $\AD:=\Aone\oplus\Atwo\oplus\Atwo'$, where $\Atwo'$ is another copy of
$\Atwo$, the multiplication is given by
    $$(\aone_1+\atwo_1+\atwo'_1)(\aone_2+\atwo_2+\atwo'_2)=\aone_1\aone_2+\atwo_1\atwo_2+(\atwo'_1\mor(\aone_2))'+(\atwo_1\atwo'_2)',$$
where $\aone\in\Aone$, $\atwo\in\Atwo$, $\atwo'\in\Atwo'$, multiplications on the r.h.s. are those of $\Aone$, $\Atwo$. This can be summarized
by saying that $\AD$ is the algebra of $2\times 2$ matrices, with elements of $\Aone$, $\Atwo$, $\Atwo'$ in the upper left, lower right, lower
left corners respectively, and zeroes in the upper right corner.

Deformation complex of $\mor:\Aone\rightarrow\Atwo$, as described in \cite{GS2} section 1, \cite{GS1}, is as follows
    $$\HC(\Aone,\Aone)\oplus\HC(\Atwo,\Atwo)\oplus\2\HC(\Aone,\Atwo),$$
where $\HC$ stands for Hochschild cochain complex, which we consider with the grading, given by identifying it with the complex of coderivations
on the bar construction. We will denote the above deformation complex by $\HC(\mor)$. It is different from $\GDefMor$, and it is not even a dg
lie algebra (in \cite{Bor} it is shown that $\HC(\mor)$ is a proper $\Li$-algebra).

In \cite{GS1} page 215, an explicit quasi-isomorphism $\tau:\HC(\mor)\rightarrow\HC(\AD)$ is constructed. We will use $\tau$ to prove that
$\Bi$-algebras $\2\BDefMor$ (definition \ref{DefinitionLie}) and $\2\HC(\AD)$ are quasi-isomorphic. First we show it when $\mor$ is injective.
\begin{proposition}\label{InjectiveCase}
    Let $\mor:\Aone\rightarrow\Atwo$ be an injective morphism of dg associative algebras, let $\AD$ be the corresponding diagram algebra.
    Then $\2\BDefMor$ and $\2\HC(\AD)$ are quasi-isomorphic $\Bi$-algebras.
\end{proposition}
\textbf{Proof:}
    We will not construct a direct quasi-isomorphism $\2\BDefMor\rightarrow\2\HC(\AD)$. Instead we will find a $\Bi$-algebra $\2\Ha$ and two
    quasi-isomorphisms of $\Bi$-algebras
    \begin{equation}\label{TwoQuis}
        \2\BDefMor\leftarrow\2\Ha\rightarrow\2\HC(\AD).
    \end{equation}
    As a dg $\K$-space $\Ha$ coincides with the dg Lie algebra $\Ha$ in the proof of theorem \ref{NonPositively}, i.e. elements of $\Ha$ are
    pairs of elements $\gOne+\gTwo\in\GOne\oplus\GTwo$, s.t.
        $$[\gaz,\gOne]+[\gaz,\gTwo]=0,$$
    where $\gaz$ is projection on $\BGG$ of the solution of the Maurer-Cartan equation in $\GMorphism$, that corresponds to $\mor$ (proposition
    \ref{MorPoint}).

    As with the dg Lie structure in the proof of theorem \ref{NonPositively}, direct computation easily shows that $\Ha$ is a
    $\Bi$-subalgebra of $\BDefMor$. Since $\mor$ is injective, inclusion $\Ha\hookrightarrow\BDefMor$ is a quasi-isomorphism. Indeed this fact
    was shown in the proof of theorem \ref{NonPositively} for inclusion of dg Lie algebras, but a $\Bi$-algebra and the corresponding dg Lie
    algebra have the same underlying complex. So the first quasi-isomorphism in formula (\ref{TwoQuis}) is constructed.

    In general it is not true that $\tau:\2\BDefMor\rightarrow\2\HC(\AD)$ is a morphism of $\Bi$-algebras. However, we will show that its
    restriction to $\Ha$ is one. First we recall the definition of $\tau$ (\cite{GS1}, page 215).

    Let $\gOne+\gTwo\in\Ha$,
    $\gOne\in\Hom((\s\Aone)\tensorM,\s\Aone)$, $\gTwo\in\Hom((\s\Atwo)^{\tensor^n},\s\Atwo)$, then
    $\tau(\gOne+\gTwo)\in\Hom(\Om(\s\AD)\tensorM,\s\AD)$ is defined as follows: for all $\aone_i\in\Aone$, $\atwo_i\in\Atwo$, $\atwo'\in\Atwo'$
        $$\tau(\gOne+\gTwo)(\aone_1,..,\aone_m):=\gOne(\aone_1,..,\aone_m),\quad
        \tau(\gOne+\gTwo)(\atwo_1,..,\atwo_n):=\gTwo(\atwo_1,..,\atwo_n),$$
        $$\tau(\gOne+\gTwo)(\atwo_1,..,\atwo_p,\atwo',\aone_1,..,\aone_q):=(\gTwo(\atwo_1,..,\atwo_p,\atwo',\mor(\aone_1),..,\mor(\aone_q)))',$$
    for all $p,q\geq 0$ s.t. $p+q+1=n$. For the rest of arguments the value of $\tau(\gOne+\gTwo)$ is set to be $0$.

    We claim that $\tau:\2\Ha\rightarrow\HC(\AD)$ is a morphism of $\Bi$-algebras. We have to show, that given $\{\gOne_i+\gTwo_i\}_{i=1}^n$ we
    have
    \begin{equation}\label{Tau1}
        \tau(\bi_{1,n-1}(\gOne_1+\gTwo_1,..,\gOne_n+\gTwo_n))=\bi_{1,n-1}(\tau(\gOne_1+\gTwo_1),..,\tau(\gOne_n+\gTwo_n)),
    \end{equation}
    \begin{equation}\label{Tau2}
        \tau(\di_n(\gOne_1+\gTwo_1,..,\gOne_n+\gTwo_n))=\di_n(\tau(\gOne_1+\gTwo_1),..,\tau(\gOne_n+\gTwo_n)),
    \end{equation}
    where we use the same symbols to denote $\Bi$-operations in both $\2\Ha$ and $\2\HC(\AD)$.

    First we prove equation (\ref{Tau1}). Since compositions of elements of $\GOne$ with elements of $\GTwo$ are always zero, the left hand side
    of this equation is
        $$\tau(\bi_{1,n-1}(\gOne_1,..,\gOne_n))+\tau(\bi_{1,n-1}(\gTwo_1,..,\gTwo_n)).$$
    It is clear that
        $$\tau(\bi_{1,n-1}(\gOne_1,..,\gOne_n))=\bi_{1,n-1}(\tau(\gOne_1),..,\tau(\gOne_n)).$$
    With the second summand it is more complicated, since $\tau(\gTwo)$ for $\gTwo\in\Hom((\s\Atwo)\tensorM,$ $\s\Atwo)$ is not zero not only on
    $(\s\Atwo)\tensorM$ but also on $(\s\Atwo)^{\tensor^{m-1}}\tensor\s\Atwo'$ and on
    $(\s\Atwo)^{\tensor^p}\tensor$ $\s\Atwo'\tensor(\s\Aone)^{\tensor^q}$ for all $p+q+1=m$. Here we use the fact that $\gOne_i+\gTwo_i\in\Ha$,
    which means that compositions of $\gOne_i$ with $\mor$ and of $\gTwo_i$ with $\mor$ are equal, i.e. $\mor_*(\gOne_i)=\mor^*(\gTwo_i)$.
    This implies that composition of  $\tau(\gOne_i)$'s with $\tau(\gTwo_1)$ is part of the image under $\tau$ of composition of $\gTwo_i$'s with
    $\gTwo_1$. So we see that equation (\ref{Tau1}) does hold.

    To prove equation (\ref{Tau2}) we first note that on $\BDefMor$ the only non-zero $\di$-operations are $\di_1$ and $\di_2$. Since by
    definition $\di_1(\Ha)\subset\Ha$, on $\Ha$ $\di_1$ coincides with the differential on $\HC(\mor)$ (\cite{GS2} section 1). We know from
    \cite{GS1} that $\tau$ maps this differential to $\di_1$ on $\HC(\AD)$, therefore it remains to show that $\di_2$ commutes with $\tau$.
    However, on $\Ha$ $\di_2$ coincides with $\bi_{1,2}(\goo+\got,-)$ and since $\goo+\got\in\Ha$, from equation (\ref{Tau1}) we conclude that
    $\di_2$ commutes with $\tau$, and hence equation (\ref{Tau2}) does hold.
$\blacksquare$

Now when we know, that for injective morphisms the $\Bi$-structure on $\2\BDefMor$ is quasi-isomorphic to $\Bi$-structure on $\2\HC(\AD)$, we
can show this for all morphisms, since these $\Bi$-structures are invariants for quasi-isomorphism classes of morphisms.
\begin{theorem}
    Let $\mor:\Aone\rightarrow\Atwo$ be a morphism of dg associative algebras. Let $\AD$ be the corresponding diagram algebra. The
    $\Bi$-algebras $\2\BDefMor$ and $\2\HC(\AD)$ are quasi-isomorphic.
\end{theorem}
\textbf{Proof:}
    Let $\mor$, $\morr$ be two quasi-isomorphic morphisms of dg associative algebras. Recall that a quasi-isomorphism between morphisms is a
    pair of quasi-isomorphisms that makes up a commutative square. Therefore we have a quasi-isomorphism between the diagram algebras $\AD$ and
    $\ADD$ of $\mor$ and $\morr$ respectively. It is well known that for quasi-isomorphic algebras, the corresponding Hochschild chain complexes
    are quasi-isomorphic $\Bi$-algebras.

    On the other hand, from lemma \ref{AlgebraChange} we know that $\2\BDefMor$ and $\2\BDefMorr$ are quasi-isomorphic. Therefore if $\2\BDefMor$
    and $\2\HC(\AD)$ are quasi-isomorphic $\Bi$-algebras, the same is true for all morphisms, that are quasi-isomorphic to $\mor$. Now the
    theorem follows from proposition \ref{InjectiveCase} since every morphism of dg associative algebras admits an injective resolution.
$\blacksquare$

In \cite{GS1} page 216 it is shown, that the Gerstenhaber algebra on the cohomology of $\HC(\AD)$ is isomorphic to the Gerstenhaber algebra on
cohomology of $\HC(\mor)=\HC(\Aone,\Aone)\oplus\HC(\Atwo,\Atwo)\oplus\HC(\Aone,\Atwo)$, that is given in \cite{GS2} section 1, page 250-251.
From the last theorem it follows that this Gerstenhaber algebra is isomorphic to the one on cohomology of $\BDefMor$.

\end{document}